\documentclass[10pt]{amsart}
\usepackage[cp1251]{inputenc}
\usepackage[english]{babel}
\usepackage{amsmath}
\usepackage{amssymb}
\usepackage{amsfonts}

\begin{document}
\renewcommand{\refname}{References}

\thispagestyle{empty}

\title[Integration Order Replacement Technique]
{Integration Order Replacement Technique for Iterated Ito Stochastic Integrals
and Iterated Stochastic Integrals With Respect to Martingales}
\author[D.F. Kuznetsov]{Dmitriy F. Kuznetsov}
\address{Dmitriy Feliksovich Kuznetsov
\newline\hphantom{iii} Peter the Great Saint-Petersburg Polytechnic University,
\newline\hphantom{iii} Polytechnicheskaya ul., 29,
\newline\hphantom{iii} 195251, Saint-Petersburg, Russia}%
\email{sde\_kuznetsov@inbox.ru}
\thanks{\sc Mathematics Subject Classification: 60H05}
\thanks{\sc Keywords: Iterated Ito stochastic integral,
Iterated stochastic integral with respect to martingales,
Integration order replacement technique, Ito formula.}

\vspace{5mm}

\maketitle {\small
\begin{quote}
\noindent{\sc Abstract.} 
The article is devoted to 
the integration order replacement technique for iterated Ito stochastic 
integrals
and iterated stochastic integrals with respect to martingales.
We consider the class of iterated Ito
stochastic integrals, for which with probability 1 the 
formulas of integration order replacement corresponding to 
the rules of classical integral calculus are correct.
The theorems on integration order replacement for 
the class of iterated Ito stochastic integrals are proved. 
Many examples of this theorems usage have been considered.
These results are generalized 
for the class of iterated stochastic integrals with respect to
martingales.
\medskip
\end{quote}
}

\vspace{4mm}

%\linespread{1.6}

\setlength{\baselineskip}{2.0em}

\tableofcontents

\setlength{\baselineskip}{1.2em}

%\linespread{1.0}

\vspace{4mm}

\section{Introduction}

\vspace{5mm}

In this article, we performed rather laborious work connected with 
the theorems on integration order replacement for  
iterated Ito stochastic integrals.
However, there may appear a
question about a practical usefulness of this theory, since the 
significant part of its conclusions directly follows from 
the Ito formula \cite{1}.

It is not difficult to see that to obtain various relations 
for iterated Ito stochastic integrals (see, for example, Sect.~6) 
using the Ito formula, first 
of all these relations should be guessed. Then it is necessary 
to introduce corresponding Ito processes and afterwards to use
the Ito formula.
It is clear that this process requires 
intellectual expenses and it is not always trivial.

On the other hand, 
the technique on integration order replacement introduced in this
article is formally comply with the similar technique for 
Riemann integrals, although it is related to Ito integrals, 
and it provides a possibility to perform transformations 
naturally (as with Riemann integrals) with iterated Ito stochastic 
integrals and to obtain 
various relations for them.

So, in order to implementation of transformations of the specific 
class of Ito processes, which is represented by iterated Ito
stochastic integrals, it is more naturally and easier to use 
the theorems on integration order replacement, than the Ito formula.

Many examples of 
these theorems usage are presented in Sect.~6. 

Note that in a lot of publications of the author \cite{77}-\cite{12aaa}
the integration order replacement technique for iterated Ito stochastic 
integrals has been successfully applied 
for the proof and development of the method of 
approximation of iterated Ito and Stratonovich stochastic integrals
based on generalized multiple Fourier series 
as well as for the construction of the so-called
unified Taylor--Ito and Taylor--Stratonovich expansions.

Let $(\Omega,{\rm F},{\sf P})$ be a complete probability
space and let $f(t,\omega):$ $[0, T]\times \Omega\rightarrow \mathbb{R}$
be the standard Wiener process
defined on the probability space $(\Omega,{\rm F},{\sf P}).$
Further, we will use the following notation:
$f(t,\omega)\stackrel{\rm def}{=}f_t$.

Let us consider the family of $\sigma$-algebras
$\left\{{\rm F}_t,\ t\in[0,T]\right\}$ defined
on the probability space $(\Omega,{\rm F},{\sf P})$ and
connected
with the Wiener process $f_t$ in such a way that

\vspace{2mm}

1.\ ${\rm F}_s\subset {\rm F}_t\subset {\rm F}$\ for
$s<t.$

\vspace{2mm}

2.\ The Wiener process $f_t$ is ${\rm F}_t$-measurable for all
$t\in[0,T].$

\vspace{2mm}

3.\ The process $f_{t+\Delta}-f_{t}$ for all
$t\ge 0,$ $\Delta>0$ is independent with
the events of $\sigma$-algebra
${\rm F}_{t}.$

\vspace{2mm}

Let us introduce the class ${\rm M}_2([0,T])$ of functions
$\xi:$ $[0,T]\times\Omega\rightarrow \mathbb{R},$ which satisfy the
conditions:

\vspace{2mm}

1. The function $\xi(t,\omega)$ is 
measurable
with respect to the pair of variables
$(t,\omega).$

\vspace{2mm}

2. The function $\xi(t,\omega)$ is ${\rm F}_t$-measurable 
for all $t\in[0,T]$ and $\xi(\tau,\omega)$ is independent 
with increments $f_{t+\Delta}-f_{t}$ 
for $t\ge \tau,$ $\Delta>0.$

\vspace{2mm}

3.\ The following relation is fulfilled
$$\int\limits_0^T{\sf M}\left\{\left(\xi(t,\omega)\right)^2\right\}dt
<\infty.
$$

\vspace{2mm}

4.\  ${\sf M}\left\{\left(\xi(t,\omega)\right)^2\right\}<\infty$
for all $t\in[0,T].$

\vspace{2mm}

For any partition 
$\tau_j^{(N)},$ $j=0, 1, \ldots, N$ of
the interval $[0,T]$ such that

\begin{equation}
\label{pr}
0=\tau_0^{(N)}<\tau_1^{(N)}<\ldots <\tau_N^{(N)}=T,\ \ \ \
\max\limits_{0\le j\le N-1}\left|\tau_{j+1}^{(N)}-\tau_j^{(N)}\right|\to 0\ \
\hbox{if}\ \ N\to \infty
\end{equation}

\vspace{3mm}
\noindent
we will define the sequense of step functions 

$$
\xi^{(N)}(t,\omega)=\xi_j(\omega)\ \ \ \hbox{w. p. 1}\ \ \ 
\hbox{for}\ \ \ t\in\left[\tau_j^{(N)},\tau_{j+1}^{(N)}\right),
$$

\vspace{3mm}
\noindent
where $j=0, 1,\ldots,N-1,$\ \ $N=1, 2,\ldots.$ 
Here and further, w.~p.~1 means 
with probability 1.

Let us define the Ito stochastic integral for
$\xi(t,\omega)\in{\rm M}_2([0,T])$ as the 
following mean-square limit \cite{1}
\begin{equation}
\label{10.1}
\hbox{\vtop{\offinterlineskip\halign{
\hfil#\hfil\cr
{\rm l.i.m.}\cr
$\stackrel{}{{}_{N\to \infty}}$\cr
}} }\sum_{j=0}^{N-1}\xi^{(N)}\left(\tau_j^{(N)},\omega\right)
\left(f\left(\tau_{j+1}^{(N)},\omega\right)-
f\left(\tau_j^{(N)},\omega\right)\right)
\stackrel{\rm def}{=}\int\limits_0^T\xi_\tau df_\tau,
\end{equation}

\vspace{2mm}
\noindent
where $\xi^{(N)}(t,\omega)$ is any step function, 
which converges
to the function $\xi(t,\omega)$
in the following sense

\vspace{-1mm}
\begin{equation}
\label{jjj}
\hbox{\vtop{\offinterlineskip\halign{
\hfil#\hfil\cr
{\rm lim}\cr
$\stackrel{}{{}_{N\to \infty}}$\cr
}} }\int\limits_0^T{\sf M}\left\{\left|
\xi^{(N)}(t,\omega)-\xi(t,\omega)\right|^2\right\}dt=0.
\end{equation}

\vspace{3mm}

It is well known \cite{1} that the Ito stochastic integral
exists as the limit (\ref{10.1}) and it does not depend on the selection 
of sequence 
$\xi^{(N)}(t,\omega)$. We suppose that standard
properties of the Ito stochastic integral are well known 
to the reader (see, for example, \cite{1}).

Let us define the stochastic integral for 
$\xi_{\tau}\in{\rm M}_2([0,T])$ as the 
following mean-square limit

\vspace{-1mm}
$$
\hbox{\vtop{\offinterlineskip\halign{
\hfil#\hfil\cr
{\rm l.i.m.}\cr
$\stackrel{}{{}_{N\to \infty}}$\cr
}} }\sum_{j=0}^{N-1}\xi^{(N)}\left(\tau_j^{(N)},\omega\right)
\left(\tau_{j+1}^{(N)}-
\tau_j^{(N)}\right)\stackrel{\rm def}{=}\int\limits_{0}^T \xi_\tau d\tau,
$$

\vspace{2mm}
\noindent
where $\xi^{(N)}(t,\omega)$ is any step function
from the class ${\rm M}_2([0,T])$,
which converges
in the sense (\ref{jjj}) 
to the function $\xi(t,\omega).$

We will introduce the class ${\rm S}_2([0,T])$ of functions 
$\xi:$ $[0,T]\times\Omega\rightarrow
\mathbb{R},$ which satisfy the conditions:

\vspace{2mm}

1.~$\xi(\tau,\omega) \in {\rm M}_2([0,T])$.

\vspace{2mm}

2.~$\xi(\tau,\omega)$
is the mean-square continuous random process at the interval
$[0,T].$

\vspace{2mm}

As we noted above,
the Ito stochastic integral exists
in the mean-square sense (see (\ref{10.1})), if the random process
$\xi(\tau,\omega)\in {\rm M}_2([0,T]),$ 
i.e., perhaps this process does not satisfy 
the property of the mean-square continuity on the interval 
$[0,T].$ In this article we will formulate and prove
the theorems on integration order replacement for the special 
class of iterated Ito stochastic integrals.
At the same time, the condition of the mean-square continuity 
of integrand in the innermost
stochastic integral will be significant.

Let us introduce the following class
of iterated stochastic integrals

\vspace{-1mm}
$$
J[\phi,\psi^{(k)}]_{T,t}=\int\limits_{t}^{T}\psi_1(t_1)\ldots
\int\limits_t^{t_{k-1}}\psi_k(t_k)\int\limits_t^{t_k}
\phi_{\tau}dw^{(k+1)}_{\tau}dw^{(k)}
_{t_k}
\ldots dw^{(1)}_{t_1},
$$

\vspace{2mm}
\noindent
where $\phi(\tau,\omega)\stackrel{\rm def}{=}\phi_{\tau},$ 
$\phi_\tau\in{\rm S}_2([t,T]),$ every
$\psi_l(\tau)$ $(l=1,\ldots,k)$ is a continous nonrandom function
at the interval $[t, T]$,
here and further 
$w_\tau^{(l)}=f_\tau$ or $w_\tau^{(l)}=\tau$
for $\tau\in[t,T]$
$(l=1,\ldots,k+1),$
$(\psi_1,\ldots,\psi_k)\stackrel{\rm def}{=}\psi^{(k)},$
$\psi^{(1)}\stackrel{\rm def}{=}\psi_1.$

We will call the stochastic integral $J[\phi,\psi^{(k)}]_{T,t}$
as the iterated Ito stochastic integral.

It is well known that for the iterated Riemann integral 
in the case of specific conditions the formula on
integration order replacement is correct.
In particular, if the nonrandom functions
$f(x)$ and $g(x)$ are continuous at the interval $[a, b],$ then

\vspace{-1mm}
\begin{equation}
\label{1.3000000}
\int\limits_a^b f(x)\int\limits_a^x g(y)dydx
=\int\limits_a^b g(y)\int\limits_y^b f(x)dxdy.
\end{equation}

\vspace{2mm}

If we suppose that for the Ito stochastic integral 

\vspace{-1mm}
$$
J[\phi,\psi_1]_{T,t}=\int\limits_{t}^{T}\psi_1(s)
\int\limits_t^{s}\phi_{\tau}dw^{(2)}_{\tau}dw^{(1)}_{s}
$$ 

\vspace{2mm}
\noindent
the formula on integration order replacement,   
which is similar to (\ref{1.3000000}), is valid, then we will have 

\vspace{-1mm}
\begin{equation}
\label{1.4000000}
\int\limits_{t}^{T}\psi_1(s)
\int\limits_t^{s}\phi_{\tau}dw^{(2)}_{\tau}dw^{(1)}
_{s}=\int\limits_{t}^T \phi_{\tau}\int\limits_{\tau}^T\psi_1(s)
dw_s^{(1)}dw_{\tau}^{(2)}.
\end{equation}

\vspace{2mm}

If, in addition 
$w_s^{(1)},\ w_s^{(2)}=f_s$ $(s\in[t, T])$ in (\ref{1.4000000}), then
the stochastic process 

\vspace{-1mm}
$$
\eta_\tau=
\phi_{\tau}\int\limits_{\tau}^T\psi_1(s)dw_s^{(1)}
$$ 

\vspace{2mm}
\noindent
does not belong to the class ${\rm M}_2([t,T]),$
and, consequently, for the Ito stochastic integral   

\vspace{-1mm}
$$
\int\limits_t^T\eta_\tau dw_\tau^{(2)}
$$

\vspace{2mm}
\noindent
on the right-hand side of (\ref{1.4000000})
the conditions of its existence are not fulfilled.

At the same time

\vspace{-1mm}
\begin{equation}
\label{rrr111}
\int\limits_t^T df_s\int\limits_t^T ds=
\int\limits_t^T (s-t)df_s+\int\limits_t^T (f_s-f_t)ds\ \ \ \hbox{w.~p.~1},
\end{equation}

\vspace{2mm}
\noindent
and we can obtain this equality, for example, using the Ito formula,
but (\ref{rrr111}) can be considered as a result of 
integration order replacement (see below). 

Actually, we can demonstrate that 

\vspace{-1mm}
$$
\int\limits_t^T (f_s-f_t)ds=
\int\limits_t^T\int\limits_t^s df_\tau ds=
\int\limits_t^T\int\limits_{\tau}^T ds df_{\tau}\ \ \ \hbox{w.~p.~1}.
$$

\vspace{2mm}

Then

\vspace{-1mm}
$$
\int\limits_t^T (s-t)df_s+\int\limits_t^T (f_s-f_t)ds=
\int\limits_t^T\int\limits_t^\tau ds df_\tau+
\int\limits_t^T\int\limits_{\tau}^T ds df_{\tau}=
\int\limits_t^T df_s\int\limits_t^T ds\ \ \ \hbox{w.~p.~1}.
$$

\vspace{3mm}

The aim of this article is to establish the strict mathematical sense 
of the formula (\ref{1.4000000}) for the case
$w_s^{(1)},$ $w_s^{(2)}=f_s$ $(s\in[t, T])$ as well as its analogue
corresponding to the iterated Ito stochastic integral 
$J[\phi,\psi^{(k)}]_{T,t},$ $k \ge 2.$
At that, we will use the definition of the Ito stochastic integral 
which is more general than (\ref{10.1}).

Let us consider the partition 
$\tau_j^{(N)},$ $j=0, 1, \ldots, N$ of
the interval $[t,T]$ such that

\vspace{-2mm}
\begin{equation}
\label{prgaba}
t=\tau_0^{(N)}<\tau_1^{(N)}<\ldots <\tau_N^{(N)}=T,\ \ \ \
\max\limits_{0\le j\le N-1}\left|\tau_{j+1}^{(N)}-\tau_j^{(N)}\right|\to 0\ \
\hbox{if}\ \ N\to \infty.
\end{equation}

\vspace{2mm}

In \cite{2} Stratonovich R.L. introduced the definition 
of the so-called 
combined stochastic integral for the specific class of 
integrated processes.
Taking this definition as a foundation, let us consider the 
following construction of stochastic integral 

\vspace{-1mm}
\begin{equation}
\label{1.5000000}
\hbox{\vtop{\offinterlineskip\halign{
\hfil#\hfil\cr
{\rm l.i.m.}\cr
$\stackrel{}{{}_{N\to \infty}}$\cr
}} }
\sum^{N-1}_{j=0} \phi_{\tau_{j}}\left(
f_{\tau_{j+1}} - f_{\tau_{j}}
\right)\theta_{\tau_{j+1}} \stackrel {{\rm def}}{=}
\int\limits_{t}^{T}\phi_{\tau}df_{\tau}\theta_{\tau},
\end{equation}

\vspace{2mm}
\noindent
where $\phi_{\tau},$ $\theta_\tau\in{\rm S}_2([t,T]),$
$\{\tau_j\}_{j=0}^{N}$ is the partition 
of the interval $[t, T],$ which satisfies to the condition (\ref{pr})
(for simplicity we write here and sometimes further 
$\tau_j$ instead of $\tau_j^{(N)}$).

Further, we will prove existence of the integral
(\ref{1.5000000}) for 
$\phi_{\tau}\in{\rm S}_2([t,T])$
and $\theta_\tau$ from a little bit narrower class of processes
than ${\rm S}_2([t,T]).$ In addition, the integral defined by
(\ref{1.5000000})
will be used for the formulation and proof 
of the theorem on integration order replacement for the iterated Ito
stochastic integrals $J[\phi,\psi^{(k)}]_{T,t},$ $k \ge 1.$

Note that under the appropriate conditions the following 
properties of stochastic integrals
defined by the formula (\ref{1.5000000}) can be proved

\vspace{-1mm}
$$
\int\limits_{t}^{T} \phi_\tau df_\tau g(\tau)=
\int\limits_{t}^{T} \phi_\tau g(\tau) df_\tau\ \ \ \hbox{w. p. 1}, 
$$

\vspace{2mm}
\noindent
where $g(\tau)$ is a continuous nonrandom function at the
interval $[t,T]$,

\vspace{-1mm}
$$
\int\limits_{t}^{T} 
\left(\alpha\phi_{\tau}+\beta\psi_{\tau}\right)df_\tau\theta_\tau=
\alpha\int\limits_{t}^{T} \phi_\tau df_\tau\theta_\tau+
\beta\int\limits_{t}^{T} \psi_{\tau} df_\tau\theta_\tau\ \ \ \hbox{w. p. 1},
$$

$$
\int\limits_{t}^{T} \phi_\tau df_\tau
\left(\alpha\theta_\tau+\beta\psi_\tau\right)=
\alpha\int\limits_{t}^{T} \phi_\tau df_\tau\theta_\tau+
\beta\int\limits_{t}^{T} \phi_\tau df_\tau\psi_\tau\ \ \ \hbox{w. p. 1},
$$

\vspace{2mm}
\noindent
where $\alpha,$ $\beta\in\mathbb{R}.$

At that, we suppose that the stochastic processes 
$\phi_\tau,$ $\theta_\tau$, and $\psi_\tau$ are such that
the integrals included in the mentioned 
properties exist.

\vspace{5mm}

\section{Formulation of the Theorem on
Integration Order Replacement
for Iterated Ito Stochastic Integrals of Arbitrary Multiplicity}

\vspace{5mm}

Let us define the stochastic integrals 
$\hat I[\psi^{(k)}]_{T,s},$ $k\ge 1$ of the form

\vspace{1mm}
$$
\hat I[\psi^{(k)}]_{T,s}=\int\limits_s^T\psi_k(t_k)dw_{t_k}^{(k)}
\int\limits_{t_k}^T\psi_{k-1}(t_{k-1})dw_{t_{k-1}}^{(k-1)}
\ldots \int\limits_{t_{2}}^T \psi_1(t_1)dw_{t_1}^{(1)}
$$

\vspace{4mm}
\noindent
in accordance with the definition (\ref{1.5000000}) by the following 
recurrence relation

\vspace{1mm}
\begin{equation}
\label{2.1000000}
\hat I[\psi^{(k)}]_{T,t}
\stackrel{\rm def}{=}
\hbox{\vtop{\offinterlineskip\halign{
\hfil#\hfil\cr
{\rm l.i.m.}\cr
$\stackrel{}{{}_{N\to \infty}}$\cr
}} }
\sum^{N-1}_{l=0} \psi_k(\tau_{l})\Delta w_{\tau_l}^{(k)} 
\hat I[\psi^{(k-1)}]_{T,\tau_{l+1}},
\end{equation}

\vspace{4mm}
\noindent
where $k\ge 1,$ 
$\hat I[\psi^{(0)}]_{T,s}\stackrel{\rm def}{=}1,$
$[s,T]\subseteq[t,T],$ here and further $\Delta w_{\tau_l}^{(i)}=
w_{\tau_{l+1}}^{(i)}-w_{\tau_l}^{(i)},$\ \
$i=1,\ldots,k+1,$\ \ $l=0, 1,\ldots,N-1.$

Then, we will define the iterated stochastic 
integral $\hat J[\phi,\psi^{(k)}]_{T,t},$\ $k\ge 1$

\vspace{1mm}
$$
\hat J[\phi,\psi^{(k)}]_{T,t}=\int\limits_{t}^T \phi_s dw_s^{(k+1)}
\hat I[\psi^{(k)}]_{T,s}
$$

\vspace{4mm}
\noindent
similarly in accordance with the definition (\ref{1.5000000})

\vspace{1mm}
$$
\hat J[\phi,\psi^{(k)}]_{T,t}
\stackrel{\rm def}{=}
\hbox{\vtop{\offinterlineskip\halign{
\hfil#\hfil\cr
{\rm l.i.m.}\cr
$\stackrel{}{{}_{N\to \infty}}$\cr
}} }
\sum^{N-1}_{l=0} \phi_{\tau_{l}}\Delta w_{\tau_l}^{(k+1)} 
\hat I[\psi^{(k)}]_{T,\tau_{l+1}}.
$$

\vspace{4mm}

Let us formulate the theorem on
integration order replacement for
iterated Ito stochastic integrals. 

\vspace{2mm} 

{\bf Theorem 1}\ \cite{3}, \cite{3aaa} 
(also see \cite{77}-\cite{20}, \cite{12a}-\cite{12aaa}). {\it Suppose that 
$\phi_\tau\in{\rm S}_2([t,T])$ and every $\psi_l(\tau)$
$(l=1,\ldots,k)$ is a continuous nonrandom function
at the interval
$[t,T]$. Then, the stochastic integral
$\hat J[\phi,\psi^{(k)}]_{T,t}$ $(k\ge 1)$ exists and 

\vspace{-1mm}
$$
J[\phi,\psi^{(k)}]_{T,t}=\hat J[\phi,\psi^{(k)}]_{T,t}\ \ \ 
\hbox{w. p. {\rm 1.}}
$$
}

\vspace{3mm}

\section{Proof of Theorem 1 for the Case of Iterated Ito Stochastic Integrals
of Multiplicity 2}

\vspace{5mm}

First, let us prove Theorem 1 for the case $k=1.$ We have 

\vspace{1mm}
$$
J[\phi,\psi_1]_{T,t}\stackrel{\rm def}{=}
\hbox{\vtop{\offinterlineskip\halign{
\hfil#\hfil\cr
{\rm l.i.m.}\cr 
$\stackrel{}{{}_{N\to \infty}}$\cr 
}} }\sum_{l=0}^{N-1} \psi_1(\tau_l)\Delta w_{\tau_l}^{(1)}\int\limits_{t}
^{\tau_l}\phi_{\tau}dw_{\tau}^{(2)}=
$$

\vspace{1mm}
\begin{equation}
\label{2.2000000}
=
\hbox{\vtop{\offinterlineskip\halign{
\hfil#\hfil\cr
{\rm l.i.m.}\cr
$\stackrel{}{{}_{N\to \infty}}$\cr
}} }\sum_{l=0}^{N-1}\psi_1(\tau_l)\Delta w_{\tau_l}^{(1)}
\sum_{j=0}^{l-1}
\int\limits_{\tau_j}
^{\tau_{j+1}}\phi_{\tau}dw_{\tau}^{(2)},
\end{equation}

\vspace{3mm}
$$
\hat J[\phi,\psi_1]_{T,t}\stackrel{\rm def}{=}
\hbox{\vtop{\offinterlineskip\halign{
\hfil#\hfil\cr
{\rm l.i.m.}\cr
$\stackrel{}{{}_{N\to \infty}}$\cr
}} }\sum_{j=0}^{N-1} \phi_{\tau_j}\Delta w_{\tau_j}^{(2)}
\int\limits_{\tau_{j+1}}^T \psi_1(s)dw_{s}^{(1)}=
$$

\vspace{1mm}
$$
=
\hbox{\vtop{\offinterlineskip\halign{
\hfil#\hfil\cr
{\rm l.i.m.}\cr
$\stackrel{}{{}_{N\to \infty}}$\cr
}} }\sum_{j=0}^{N-1}\phi_{\tau_j}\Delta w_{\tau_j}^{(2)}
\sum_{l=j+1}^{N-1}\int\limits_{\tau_{l}}^{\tau_{l+1}}
\psi_1(s)dw_{s}^{(1)}=
$$

\vspace{1mm}
\begin{equation}
\label{2.3000000}
=\hbox{\vtop{\offinterlineskip\halign{
\hfil#\hfil\cr
{\rm l.i.m.}\cr
$\stackrel{}{{}_{N\to \infty}}$\cr
}} }\sum_{l=0}^{N-1} \int\limits_{\tau_{l}}^{\tau_{l+1}}
\psi_1(s)dw_{s}^{(1)}
\sum_{j=0}^{l-1}
\phi_{\tau_j}\Delta w_{\tau_j}^{(2)}.
\end{equation}

\vspace{5mm}

It is clear that if the difference $\varepsilon_N$ of
prelimit expressions  
on the right-hand sides of (\ref{2.2000000}) and  
(\ref{2.3000000}) tends
to zero when $N\to\infty$ in the mean-square sense, 
then the stochastic integral  
$\hat J[\phi,\psi_1]_{T,t}$ exists 
and 

\vspace{-2mm}
$$
J[\phi,\psi_1]_{T,t}=\hat J[\phi,\psi_1]_{T,t}\ \ \ \hbox{\rm w. p. 1.}
$$

\vspace{5mm}

The difference $\varepsilon_N$ can be presented in the form 
$\varepsilon_N=\tilde\varepsilon_N+\hat\varepsilon_N$, where

\vspace{2mm}
$$
\tilde\varepsilon_N=\sum_{l=0}^{N-1} \psi_1(\tau_l)\Delta w_{\tau_l}^{(1)}
\sum_{j=0}^{l-1}
\int\limits_{\tau_j}
^{\tau_{j+1}}\left(\phi_{\tau}-\phi_{\tau_j}\right)dw_{\tau}^{(2)};
$$

\vspace{2mm}
$$
\hat\varepsilon_N=
\sum_{l=0}^{N-1} \int\limits_{\tau_{l}}^{\tau_{l+1}}
\left(\psi_1(\tau_l)-\psi_1(s)\right)dw_{s}^{(1)}
\sum_{j=0}^{l-1}
\phi_{\tau_j}\Delta w_{\tau_j}^{(2)}.
$$

\vspace{5mm}

We will demonstrate that  

\vspace{-3mm}
$$
\hbox{\vtop{\offinterlineskip\halign{
\hfil#\hfil\cr
{\rm l.i.m.}\cr
$\stackrel{}{{}_{N\to \infty}}$\cr
}} }\varepsilon_N=0.
$$ 

\vspace{5mm}

In order to do it we will analyze four cases:

\vspace{2mm}

1.\ $w_{\tau}^{(2)}=f_{\tau},\ \Delta w_{\tau_l}^{(1)}=\Delta f_{\tau_l}.$

\vspace{2mm}

2.\ $w_{\tau}^{(2)}=\tau,\ \Delta w_{\tau_l}^{(1)}=\Delta f_{\tau_l}.$

\vspace{2mm}

3.\ $w_{\tau}^{(2)}=f_{\tau},\ \Delta w_{\tau_l}^{(1)}=\Delta \tau_l.$

\vspace{2mm}

4.\ $w_{\tau}^{(2)}=\tau,\ \Delta w_{\tau_l}^{(1)}=\Delta \tau_l.$

\vspace{4mm}

Consider the well known standard moment properties 
of stochastic integrals \cite{1}

$$
{\sf M}\left\{\left|\int\limits_{t_0}^t \xi_\tau
df_\tau\right|^{2}\right\}=
\int\limits_{t_0}^t {\sf M}\left\{|\xi_\tau|^{2}\right\}d\tau,
$$

\begin{equation}
\label{99.02}
{\sf M}\left\{\left|\int\limits_{t_0}^t \xi_\tau
d\tau\right|^{2}\right\} \le (t-t_0)
\int\limits_{t_0}^t {\sf M}\left\{|\xi_\tau|^{2}\right\}d\tau,
\end{equation}

\vspace{2mm}
\noindent
where 
$\xi_{\tau}\in{\rm M}_2([t_0,t]).$

For Case 1 using 
standard moment properties for the Ito
stochastic integral as well as
mean-square continuity (which means
uniform mean-square continuity) 
of the process $\phi_\tau$ on the interval $[t, T]$, we obtain

$$
{\sf M}\left\{\left|\tilde\varepsilon_N\right|^2\right\}=
\sum_{k=0}^{N-1} \psi_1^2(\tau_k)\Delta \tau_k
\sum_{j=0}^{k-1}
\int\limits_{\tau_j}
^{\tau_{j+1}}{\sf M}\left\{\left|\phi_{\tau}-\phi_{\tau_j}\right|^2
\right\}d{\tau}<
$$

\vspace{2mm}
$$
<C^2\varepsilon
\sum_{k=0}^{N-1}\Delta \tau_k \sum_{j=0}^{k-1}\Delta \tau_j
<C^2\varepsilon\frac{(T-t)^2}{2},
$$

\vspace{4mm}
\noindent 
i.e. ${\sf M}\left\{\left|\tilde\varepsilon_N\right|^2\right\}
\to 0$ when
$N\to \infty.$ Here $\Delta\tau_j<\delta(\varepsilon),$
$j=0, 1,\ldots,N-1$
($\delta(\varepsilon)>0$ exists for any
$\varepsilon>0$ and it does not depend on $\tau$),
$|\psi_1(\tau)|<C.$

Let us consider Case 2. Using the Minkowski inequality,
uniform mean-square continuity of the process $\phi_\tau$
as well as the estimate (\ref{99.02}) for the stochastic integral,
we have

$$
{\sf M}\left\{\left|\tilde\varepsilon_N\right|^2\right\}=
\sum_{k=0}^{N-1} \psi_1^2(\tau_k)\Delta \tau_k
{\sf M}\left\{\left(\sum_{j=0}^{k-1}
\int\limits_{\tau_j}^{\tau_{j+1}}(\phi_{\tau}-\phi_{\tau_j})d\tau\right)^2
\right\}\le
$$

\vspace{2mm}
$$
\le 
\sum_{k=0}^{N-1} \psi_1^2(\tau_k)\Delta \tau_k
\left(\sum_{j=0}^{k-1}\left(
{\sf M}\left\{\left(
\int\limits_{\tau_j}^{\tau_{j+1}}(\phi_{\tau}-\phi_{\tau_j})d\tau\right)^2
\right\}\right)^{1/2}\ \right)^2<
$$

\vspace{2mm}
$$
< C^2\varepsilon
\sum_{k=0}^{N-1}\Delta \tau_k \left(\sum_{j=0}^{k-1}\Delta \tau_j\right)^2
<C^2\varepsilon \frac{(T-t)^3}{3},
$$

\vspace{4mm}
\noindent
i.e. ${\sf M}\left\{\left|\tilde\varepsilon_N\right|^2\right\}
\to 0$ 
when  
$N\to \infty.$ Here 
$\Delta\tau_j<\delta(\varepsilon),$ $j=0, 1,\ldots,N-1$
($\delta(\varepsilon)>0$ exists for any
$\varepsilon>0$ and it does not depend on $\tau$),
$|\psi_1(\tau)|< C.$

For Case 3 using the Minkowski inequality,
standard moment properties for the Ito stochastic integral
as well as uniform mean-square continuity 
of the process $\phi_\tau$, we find

$$
{\sf M}\left\{\left|\tilde\varepsilon_N\right|^2\right\}\le
\left(
\sum_{k=0}^{N-1}\left| \psi_1(\tau_k)\right| \Delta \tau_k
\left({\sf M}\left\{\left(\sum_{j=0}^{k-1}
\int\limits_{\tau_j}^{\tau_{j+1}}(\phi_{\tau}-\phi_{\tau_j})df_{\tau}\right)^2
\right\}\right)^{1/2}\ \right)^2=
$$

\vspace{2mm}
$$
=
\left(
\sum_{k=0}^{N-1}| \psi_1(\tau_k)| \Delta \tau_k
\left(\sum_{j=0}^{k-1}
\int\limits_{\tau_j}^{\tau_{j+1}}{\sf M}\left\{|\phi_{\tau}-\phi_{\tau_j}|
^2\right\}d\tau\right)^{1/2}\ \right)^2 <
$$

\vspace{2mm}
$$
< C^2\varepsilon
\left(\sum_{k=0}^{N-1}\Delta \tau_k \left(\sum_{j=0}^{k-1}\Delta \tau_j 
\right)^{1/2}\
\right)^2
<C^2\varepsilon\frac{4(T-t)^3}{9},
$$

\vspace{4mm}
\noindent
i.e. ${\sf M}\left\{\left|\tilde\varepsilon_N\right|^2\right\}
\to 0$ 
when  
$N\to \infty.$ Here 
$\Delta\tau_j<\delta(\varepsilon),$ $j=0, 1,\ldots,N-1$
($\delta(\varepsilon)>0$ exists for any
$\varepsilon>0$ and it does not depend on $\tau$),
$|\psi_1(\tau)|< C.$

Finally, for Case 4 using 
the Minkowski inequality, 
uniform mean-square continuity of the process $\phi_\tau$
as well as the estimate (\ref{99.02}) for the stochastic 
integral, 
we obtain

$$
{\sf M}\left\{\left|\tilde\varepsilon_N\right|^2\right\}\le
\left(
\sum_{k=0}^{N-1} \sum_{j=0}^{k-1}| \psi_1(\tau_k)| \Delta \tau_k
\left({\sf M}\left\{\left(
\int\limits_{\tau_j}^{\tau_{j+1}}(\phi_{\tau}-\phi_{\tau_j})d\tau\right)^2
\right\}\right)^{1/2}\ \right)^2 <
$$

\vspace{2mm}
$$
< C^2\varepsilon
\left(\sum_{k=0}^{N-1}\Delta \tau_k \sum_{j=0}^{k-1}\Delta \tau_j
\right)^2
<C^2\varepsilon\frac{(T-t)^4}{4},
$$

\vspace{4mm}
\noindent
i.e. ${\sf M}\left\{\left|\tilde\varepsilon_N\right|^2\right\} \to 0$
when  
$N\to \infty.$ Here 
$\Delta\tau_j<\delta(\varepsilon),$ $j=0, 1,\ldots,N-1$
($\delta(\varepsilon)>0$ exists for any
$\varepsilon>0$ and it does not depend on $\tau$),\
$|\psi_1(\tau)|< C.$

Thus, we have proved that  

\vspace{-3mm}
$$
\hbox{\vtop{\offinterlineskip\halign{
\hfil#\hfil\cr
{\rm l.i.m.}\cr
$\stackrel{}{{}_{N\to \infty}}$\cr
}} }\tilde\varepsilon_N=0.
$$ 

\vspace{3mm}

Analogously, taking into account 
the uniform continuity of the function $\psi_1(\tau)$
on the interval $[t, T]$, we can demonstrate that 

\vspace{-3mm}
$$
\hbox{\vtop{\offinterlineskip\halign{
\hfil#\hfil\cr
{\rm l.i.m.}\cr
$\stackrel{}{{}_{N\to \infty}}$\cr
}} }\hat\varepsilon_N=0.
$$ 

\vspace{3mm}

Consequently,

\vspace{-3mm}
$$
\hbox{\vtop{\offinterlineskip\halign{
\hfil#\hfil\cr
{\rm l.i.m.}\cr
$\stackrel{}{{}_{N\to \infty}}$\cr
}} }\varepsilon_N=0.
$$ 

\vspace{3mm}

Theorem 1 is proved for the case $k=1$.

\vspace{2mm}

{\bf Remark 1.}\ {\it Proving Theorem {\rm 1,} we used the fact that 
if the stochastic process 
$\phi_t$ is
mean-square continuous at the interval $[t, T],$ 
then it is uniformly mean-square continuous at
this interval, i.e. $\forall$ $\varepsilon>0$
$\exists$ $\delta(\varepsilon)>0$ such that
for all $t_1, t_2\in [t, T]$ satisfying 
the condition $|t_1-t_2|<\delta(\varepsilon)$ the inequality 

\vspace{-1mm}
$$
{\sf M}\left\{\left|\phi_{t_1}-\phi_{t_2}\right|^2\right\}<\varepsilon
$$

\vspace{2mm}
\noindent
is fulfilled
{\rm (}here $\delta(\varepsilon)$ does not depend on $t_1$ and $t_2${\rm )}.

\vspace{2mm}

{\bf Proof.}\ Suppose that the stochastic process $\phi_t$ is mean-square 
continuous at the interval $[t, T]$, but not 
uniformly mean-square continuous 
at this interval. Then for some $\varepsilon>0$
and $\forall$ $\delta(\varepsilon)>0$ $\exists$ $t_1, t_2\in[t, T]$ 
such that $|t_1-t_2|<\delta(\varepsilon),$ but

\vspace{-1mm}
$$
{\sf M}\left\{\left|\phi_{t_1}-\phi_{t_2}\right|^2\right\}\ge\varepsilon.
$$

\vspace{2mm}

Consequently, for $\delta=\delta_n=1/n$ $(n\in{\bf N})$
$\exists$ $t_1^{(n)},$ $t_2^{(n)}\in[t, T]$
such that 

$$
\left|t_1^{(n)}-t_2^{(n)}\right|<\frac{1}{n},
$$

\vspace{2mm}
\noindent
but 
$$
{\sf M}\left\{\left|\phi_{t_1^{(n)}}-\phi_{t_2^{(n)}}\right|^2\right\}
\ge\varepsilon.
$$

\vspace{2mm}

The sequence $t_1^{(n)}$ $(n\in{\bf N})$
is bounded, consequently, according to the Bol\-za\-no--Wei\-er\-strass 
Theorem,
we can choose from it the 
subsequence $t_1^{(k_n)}$ $(n\in{\bf N})$
that converges
to a certain number
$\tilde t$ {\rm (}it is simple to demonstrate that 
$\tilde t\in[t, T]${\rm )}.
Similarly to it and in virtue of the inequality 

\vspace{-1mm}
$$
\left|t_1^{(n)}-t_2^{(n)}\right|<\frac{1}{n}
$$

\vspace{3mm}
\noindent
we have $t_2^{(k_n)}\to\tilde t$ when $n\to\infty.$

According to the mean-square continuity of the process $\phi_t$
at the moment $\tilde t$ and the elementary inequality
$(a+b)^2\le 2(a^2+b^2)$, we obtain

\vspace{-2mm}
$$
0\le {\sf M}\left\{\left|\phi_{t_1^{(k_n)}}-\phi_{t_2^{(k_n)}}
\right|^2\right\}\le
$$

$$
\le
2\left({\sf M}\left\{\left|\phi_{t_1^{(k_n)}}-\phi_{\tilde t}\right|^2\right\}+
{\sf M}\left\{\left|\phi_{t_2^{(k_n)}}-\phi_{\tilde t}\right|^2\right\}
\right)\to 0
$$

\vspace{4mm}
\noindent
when $n\to\infty.$ Then 
$$
\lim\limits_{n\to\infty}{\sf M}
\left\{\left|\phi_{t_1^{(k_n)}}-\phi_{t_2^{(k_n)}}\right|^2\right\}=0.
$$

\vspace{4mm}

It is impossible by virtue of the fact that 

\vspace{-1mm}
$$
{\sf M}\left\{\left|\phi_{t_1^{(k_n)}}-
\phi_{t_2^{(k_n)}}\right|^2\right\}\ge\varepsilon>0.
$$

\vspace{4mm}

The obtained contradiction proves the required statement.
}

\vspace{5mm}

\section{Proof of Theorem 1 for the Case of Iterated Ito Stochastic Integrals
of Arbitrary Multiplicity}

\vspace{5mm}

Let us prove Theorem 1 for the case $k>1.$
In order to do it  
we will introduce the following notations

\vspace{1mm}
$$
I[\psi_q^{(r+1)}]_{\theta,s}
\stackrel{\rm def}{=}
\int\limits_{s}^{\theta}\psi_{q}(t_1)\ldots
\int\limits_{s}^{t_{r}} \psi_{q+r}(t_{r+1})
dw_{t_{r+1}}^{(q+r)}
\ldots dw_{t_1}^{(q)},
$$

\vspace{2mm}
$$
J[\phi,\psi_q^{(r+1)}]_{\theta,s}
\stackrel{\rm def}{=}
\int\limits_{s}^{\theta}\psi_{q}(t_1)\ldots
\int\limits_{s}^{t_{r}} \psi_{q+r}(t_{r+1})
\int\limits_s^{t_{r+1}}\phi_\tau dw_\tau^{(q+r+1)}
dw_{t_{r+1}}^{(q+r)}
\ldots dw_{t_1}^{(q)},
$$

\vspace{2mm}
$$
G[\psi_q^{(r+1)}]_{n,m}=\sum_{j_q=m}^{n-1}\sum_{j_{q+1}=m}^{j_q-1}
\ldots \sum_{j_{q+r}=m}^{j_{q+r-1}-1}
\prod_{l=q}^{r+q}I[\psi_l]_{\tau_{j_l+1},\tau_{j_l}},
$$

\vspace{4mm}
$$
(\psi_q,\ldots,\psi_{q+r})\stackrel{\rm def}{=}\psi_q^{(r+1)},\ \ \
\psi_q^{(1)}\stackrel{\rm def}{=}\psi_q,
$$

\vspace{1mm}
$$
(\psi_1,\ldots,\psi_{r+1})\stackrel{\rm def}{=}\psi_1^{(r+1)},\ \ \ 
\psi_1^{(r+1)}\stackrel{\rm def}{=}\psi^{(r+1)}
$$

\vspace{7mm}

Note that according to notations introduced above

$$
I[\psi_l]_{s,\theta}=\int\limits_{\theta}^{s}
\psi_l(\tau)dw_{\tau}^{(l)}.
$$

\vspace{3mm}

To prove Theorem 1 for $k>1$
it is enough to show
that

\vspace{1mm}
\begin{equation}
\label{2.4000000}
J[\phi,\psi^{(k)}]_{T,t}=
\hbox{\vtop{\offinterlineskip\halign{
\hfil#\hfil\cr
{\rm l.i.m.}\cr
$\stackrel{}{{}_{N\to \infty}}$\cr
}} }
S[\phi,\psi^{(k)}]_N=
\hat J[\phi,\psi^{(k)}]_{T,t}\ \ \ \hbox{w.~p.~1},
\end{equation}

\vspace{3mm}
\noindent
where
$$
S[\phi,\psi^{(k)}]_N=
G[\psi^{(k)}]_{N,0}
\sum_{l=0}^{j_k-1}\phi_{\tau_l}\Delta w_{\tau_l}^{(k+1)},
$$

\vspace{3mm}
\noindent
where $\Delta w_{\tau_l}^{(k+1)}=w_{\tau_{l+1}}^{(k+1)}-w_{\tau_l}^{(k+1)}$.

First, let us prove the right equality in (\ref{2.4000000}).
We have

\begin{equation}
\label{2.5000000}
\hat J[\phi,\psi^{(k)}]_{T,t}
\stackrel{\rm def}{=}
\hbox{\vtop{\offinterlineskip\halign{
\hfil#\hfil\cr
{\rm l.i.m.}\cr
$\stackrel{}{{}_{N\to \infty}}$\cr
}} }
\sum^{N-1}_{l=0} \phi_{\tau_l}\Delta w_{\tau_l}^{(k+1)}
\hat I[\psi^{(k)}]_{T,\tau_{l+1}}.
\end{equation}

\vspace{4mm}

On the basis of the inductive hypothesis we obtain that

\begin{equation}
\label{2.6000000}
I[\psi^{(k)}]_{T,\tau_{l+1}}=\hat I[\psi^{(k)}]_{T,\tau_{l+1}}\ \ \ 
\hbox{w. p. 1},
\end{equation}

\vspace{4mm}
\noindent
where $\hat I[\psi^{(k)}]_{T,s}$ 
is defined in accordance with (\ref{2.1000000}) and

\vspace{1mm}
$$
I[\psi^{(k)}]_{T,s}=\int\limits_{s}^T \psi_1(t_1)\ldots
\int\limits_{s}^{t_{k-2}} \psi_{k-1}(t_{k-1})
\int\limits_{s}^{t_{k-1}} \psi_k(t_k)dw_{t_k}^{(k)}dw_{t_{k-1}}^{(k-1)}
\ldots dw_{t_1}^{(1)}.
$$

\vspace{4mm}

Let us note that when $k\ge 4$ (for $k=2,$ $3$ the arguments are similar)
due to additivity of the Ito stochastic 
integral the following equalities are correct

\vspace{1mm}
$$
I[\psi^{(k)}]_{T,\tau_{l+1}}
=\sum_{j_1=l+1}^{N-1}\int\limits_{\tau_{j_1}}^{\tau_{j_1+1}}
\psi_1(t_1)\int\limits_{\tau_{l+1}}^{t_1}\psi_{2}(
t_2)I[\psi_3^{(k-2)}]_{t_2,\tau_{l+1}}
dw_{t_2}^{(2)}dw_{t_1}^{(1)}=
$$

\vspace{3mm}
$$
=\sum_{j_1=l+1}^{N-1}\int\limits_{\tau_{j_1}}^{\tau_{j_1+1}}
\psi_1(t_1)\left(\sum_{j_2=l+1}^{j_1-1}
\int\limits_{\tau_{j_2}}^{\tau_{j_2+1}}
+\int\limits_{\tau_{j_1}}^{t_1}\right) 
\psi_{2}(t_2)I[\psi_3^{(k-2)}]_{t_2,\tau_{l+1}}
dw_{t_2}^{(2)}dw_{t_1}^{(1)}=
$$

\vspace{3mm}
\begin{equation}
\label{2.7000000}
=\ldots=G[\psi^{(k)}]_{N,l+1}+H[\psi^{(k)}]_{N,l+1}\ \ \ \hbox{w.~p.~1},
\end{equation}

\vspace{5mm}
\noindent
where

\vspace{1mm}
$$
H[\psi^{(k)}]_{N,l+1}=
\sum_{j_1=l+1}^{N-1}\int\limits_{\tau_{j_1}}^{\tau_{j_1+1}}
\psi_1(s)\int\limits_{\tau_{j_1}}^{s}
\psi_2(\tau)I[\psi_3^{(k-2)}]_{\tau,\tau_{l+1}}dw_\tau^{(2)}
dw_s^{(1)} +
$$

\vspace{3mm}
$$
+ \sum_{r=2}^{k-2}
G[\psi^{(r-1)}]_{N,l+1}
\sum_{j_r=l+1}^{j_{r-1}-1}\int\limits_{\tau_{j_r}}^{\tau_{j_r+1}}
\psi_r(s)\int\limits_{\tau_{j_r}}^{s}
\psi_{r+1}(\tau)I[\psi_{r+2}^{(k-r-1)}]_{\tau,\tau_{l+1}}
dw_\tau^{(r+1)}
dw_s^{(r)} +
$$

\vspace{2mm}
\begin{equation}
\label{2.8000000}
+ G[\psi^{(k-2)}]_{N,l+1}
\sum_{j_{k-1}=l+1}^{j_{k-2}-1}I[\psi_{k-1}^{(2)}]_{\tau_{j_{k-1}+1},
\tau_{j_{k-1}}}.
\end{equation}

\vspace{7mm}

Let us substitute (\ref{2.7000000}) into (\ref{2.6000000}) and
(\ref{2.6000000}) into (\ref{2.5000000}). Then w.~p.~1

\begin{equation}
\label{2.9000000}
\hat J[\phi,\psi^{(k)}]_{T,t}=
\hbox{\vtop{\offinterlineskip\halign{
\hfil#\hfil\cr
{\rm l.i.m.}\cr
$\stackrel{}{{}_{N\to \infty}}$\cr
}} }
\sum_{l=0}^{N-1} \phi_{\tau_l}\Delta w_{\tau_l}^{(k+1)}
\biggl(G[\psi^{(k)}]_{N,l+1}+
H[\psi^{(k)}]_{N,l+1}\biggr).
\end{equation}

\vspace{3mm}

Since
\begin{equation}
\label{2.10000001}
\sum_{j_1=0}^{N-1} \sum_{j_2=0}^{j_1-1}\ldots
\sum_{j_k=0}^{j_{k-1}-1} 
a_{j_1\ldots j_k}=\sum_{j_k=0}^{N-1} 
\sum_{j_{k-1}=j_k+1}^{N-1}\ldots \sum_{j_1=j_2+1}^{N-1}
a_{j_1\ldots j_k},
\end{equation}

\vspace{3mm}
\noindent
where $a_{j_1\ldots j_k}$ are scalars,
then

\vspace{-2mm}
\begin{equation}
\label{2.11000000}
G[\psi^{(k)}]_{N,l+1}=
\sum^{N-1}_{j_k=l+1}\ldots
\sum_{j_1=j_2+1}^{N-1}
\prod_{l=1}^k I[\psi_l]_{\tau_{j_l+1},\tau_{j_l}}.
\end{equation}

\vspace{3mm}

Let us substitute (\ref{2.11000000}) into 
$$
\sum\limits_{l=0}^{N-1}
\phi_{\tau_l}\Delta w_{\tau_l}^{(k+1)}
G[\psi^{(k)}]_{N,l+1}
$$

\vspace{4mm}
\noindent
and again use 
the formula (\ref{2.10000001}). Then 

\vspace{-1mm}
\begin{equation}
\label{2.12000000}
\sum\limits_{l=0}^{N-1}
\phi_{\tau_l}\Delta w_{\tau_l}^{(k+1)}
G[\psi^{(k)}]_{N,l+1}=S[\phi,\psi^{(k)}]_{N}.
\end{equation}

\vspace{4mm}

Let us suppose that the limit
\begin{equation}
\label{li}
\hbox{\vtop{\offinterlineskip\halign{
\hfil#\hfil\cr
{\rm l.i.m.}\cr
$\stackrel{}{{}_{N\to \infty}}$\cr
}} }S[\phi,\psi^{(k)}]_{N}
\end{equation} 

\vspace{4mm}
\noindent
exists (its existence will be proved further).

Then from (\ref{2.12000000}) and (\ref{2.9000000}) it follows that
for proof of the right equality in  
(\ref{2.4000000}) we have to demonstrate that w.~p.~1
\begin{equation}
\label{2.13000000}
\hbox{\vtop{\offinterlineskip\halign{
\hfil#\hfil\cr
{\rm l.i.m.}\cr
$\stackrel{}{{}_{N\to \infty}}$\cr
}} }
\sum_{l=0}^{N-1}
\phi_{\tau_l}\Delta w_{\tau_l}^{(k+1)}
H[\psi^{(k)}]_{N,l+1}=0.
\end{equation}

\vspace{4mm}

Analyzing the second moment of the prelimit expression on the 
left-hand side of (\ref{2.13000000}) and taking into account
(\ref{2.8000000}), the independence of
$\phi_{\tau_l},$ $\Delta w_{\tau_l}^{(k+1)},$ and
$H[\psi^{(k)}]_{N,l+1}$ as well as the
standard estimates for second moments
of stochastic integrals and the Minkowski inequality,
we find that (\ref{2.13000000})
is correct.
Thus, by
the assumption of existence of the limit 
(\ref{li}) 
we obtain that the right equality in (\ref{2.4000000}) is fulfilled.

Let us demonstrate that the left equality in  
(\ref{2.4000000}) is also fulfilled.

We have 
\begin{equation}
\label{2.14000000}
J[\phi,\psi^{(k)}]_{T,t}\stackrel{\rm def}{=}
\hbox{\vtop{\offinterlineskip\halign{
\hfil#\hfil\cr
{\rm l.i.m.}\cr
$\stackrel{}{{}_{N\to \infty}}$\cr
}} }\sum_{l=0}^{N-1}\psi_1(\tau_l)\Delta w_{\tau_l}^{(1)} 
J[\phi,\psi_2^{(k-1)}]_{\tau_l,t}.
\end{equation}

\vspace{2mm}

Let us use for
the integral $J[\phi,\psi_2^{(k-1)}]_{\tau_l,t}$
in (\ref{2.14000000}) the same arguments, 
which resulted to the relation
(\ref{2.7000000})
for the integral
$I[\psi^{(k)}]_{T,\tau_{l+1}}$. 
After that let us substitute 
the expression obtained for
the integral $J[\phi,\psi_2^{(k-1)}]_{\tau_l,t}$
into (\ref{2.14000000}). 

Further, using 
the Minkowski inequality and standard 
estimates for second moments of stochastic integrals
it is easy to obtain that

\vspace{-1mm}
\begin{equation}
\label{2.15000000}
J[\phi,\psi^{(k)}]_{T,t}=
\hbox{\vtop{\offinterlineskip\halign{
\hfil#\hfil\cr
{\rm l.i.m.}\cr
$\stackrel{}{{}_{N\to \infty}}$\cr
}} }R[\phi,\psi^{(k)}]_N\ \ \ \hbox{w. p. 1},
\end{equation}

\vspace{2mm}
\noindent
where
$$
R[\phi,\psi^{(k)}]_N=
\sum^{N-1}_{j_1=0}
\psi_1(\tau_{j_1})\Delta w_{\tau_{j_1}}^{(1)}
G[\psi_2^{(k-1)}]_{j_1,0}
\sum_{l=0}^{j_k-1}
\int\limits_{\tau_l}^{\tau_{l+1}}\phi_{\tau}dw_{\tau}^{(k+1)}.
$$

\vspace{2mm}

We will demonstrate that 

\vspace{-1mm}
\begin{equation}
\label{2.16000000}
\hbox{\vtop{\offinterlineskip\halign{
\hfil#\hfil\cr
{\rm l.i.m.}\cr
$\stackrel{}{{}_{N\to \infty}}$\cr
}} }R[\phi,\psi^{(k)}]_N=
\hbox{\vtop{\offinterlineskip\halign{
\hfil#\hfil\cr
{\rm l.i.m.}\cr
$\stackrel{}{{}_{N\to \infty}}$\cr
}} }S[\phi,\psi^{(k)}]_N\ \ \ \hbox{w. p. 1}.
\end{equation}

\vspace{3mm}

It is easy to see that  

\begin{equation}
\label{2.17000000}
R[\phi,\psi^{(k)}]_N=
U[\phi,\psi^{(k)}]_N+
V[\phi,\psi^{(k)}]_N+
S[\phi,\psi^{(k)}]_N\ \ \ \hbox{w. p. 1},
\end{equation}

\vspace{4mm}
\noindent
where

\vspace{-2mm}
$$
U[\phi,\psi^{(k)}]_N=
\sum^{N-1}_{j_1=0}
\psi_1(\tau_{j_1})\Delta w_{\tau_{j_1}}^{(1)}
G[\psi_2^{(k-1)}]_{j_1,0}
\sum_{l=0}^{j_k-1}
I[\Delta\phi]_{\tau_{l+1},\tau_l},
$$

\vspace{3mm}
$$
V[\phi,\psi^{(k)}]_N=
\sum^{N-1}_{j_1=0}
I[\Delta\psi_1]_{\tau_{j_1+1},\tau_{j_1}}
G[\psi_2^{(k-1)}]_{j_1,0}
\sum_{l=0}^{j_k-1}
\phi_{\tau_l}\Delta w_{\tau_l}^{(k+1)},
$$

\vspace{3mm}
$$
I[\Delta\psi_1]_{\tau_{j_1+1},\tau_{j_1}}
=\int\limits_{\tau_{j_1}}^{\tau_{j_1+1}}
(\psi_1(\tau_{j_1})-\psi_1(\tau))dw_{\tau}^{(1)},
$$

\vspace{2mm}
$$
I[\Delta\phi]_{\tau_{l+1},\tau_l}=\int\limits_{\tau_l}^{\tau_{l+1}}
(\phi_{\tau}-\phi_{\tau_l})dw_{\tau}^{(k+1)}.
$$

\vspace{4mm}

Using 
the Minkowski inequality, standard 
estimates for second moments of stochastic integrals,
the condition that the process $\phi_\tau$ belongs
to the class ${\rm S}_2([t,T])$
as well as   
continuity (which means uniform continuity)
of the function $\psi_1(\tau)$,
we obtain that

\vspace{1mm}
$$
\hbox{\vtop{\offinterlineskip\halign{
\hfil#\hfil\cr
{\rm l.i.m.}\cr
$\stackrel{}{{}_{N\to \infty}}$\cr
}} }V[\phi,\psi^{(k)}]_N=
\hbox{\vtop{\offinterlineskip\halign{
\hfil#\hfil\cr
{\rm l.i.m.}\cr
$\stackrel{}{{}_{N\to \infty}}$\cr
}} }U[\phi,\psi^{(k)}]_N=0\ \ \ \hbox{w. p. 1}.
$$

\vspace{4mm}

Then, considering (\ref{2.17000000}), we obtain
(\ref{2.16000000}). From (\ref{2.16000000}) and (\ref{2.15000000}) 
it follows that 
the left equality in (\ref{2.4000000}) is fulfilled.

Note that the limit  
(\ref{li})
exists because it is equal to the stochastic
integral $J[\phi,\psi^{(k)}]_{T,t}$,
which exists under the conditions of Theorem 1.
So, the chain of equalities (\ref{2.4000000}) is proved. 
Theorem 1 is proved.

\vspace{5mm}

\section{Corollaries and Generalizations of Theorem 1}

\vspace{5mm}

Denote $D_k=\{(t_1,\ldots,t_k): t\le t_1<\ldots<t_k\le T\}$.
We will use the same symbol $D_{k}$ to denote the open and closed 
domains corresponding to the domain $D_{k}$.
However, we always specify what domain we consider (open or closed).

Suppose that the following conditions are fulfilled:

\vspace{2mm}

AI.\ $\xi_{\tau}\in{\rm S}_2([t,T]).$

\vspace{2mm}

AII.\ $\Phi(t_1,\ldots,{t_{k-1}})$ is a continuous nonrandom function 
in the closed domain $D_{k-1}$. 

\vspace{2mm}

Let us define the following stochastic integrals

\vspace{1mm}
$$
\hat J[\xi,\Phi]_{T,t}^{(k)}=\int\limits_t^T \xi_{t_k}d{\bf w}_{t_k}^{(i_k)}
\ldots 
\int\limits_{t_{3}}^T d{\bf w}_{t_{2}}^{(i_{2})}
\int\limits_{t_{2}}^{T}\Phi(t_1,t_2,\ldots,t_{k-1})
d{\bf w}_{t_1}^{(i_1)}
\stackrel{\rm def}{=}
$$

\vspace{2mm}
$$
\stackrel{\rm def}{=}\hbox{\vtop{\offinterlineskip\halign{
\hfil#\hfil\cr
{\rm l.i.m.}\cr
$\stackrel{}{{}_{N\to \infty}}$\cr
}} }\sum_{l=0}^{N-1}\xi_{\tau_l}\Delta{\bf w}_{\tau_{l}}^{(i_k)}
\int\limits_{\tau_{l+1}}^T
d{\bf w}_{t_{k-1}}^{(i_{k-1})}
\ldots\int\limits_{t_{3}}^T
d{\bf w}_{t_{2}}^{(i_{2})}
\int\limits_{t_{2}}^T
\Phi(t_1,t_2,\ldots,t_{k-1})d{\bf w}_{t_1}^{(i_1)}
$$

\vspace{4mm}
\noindent
for $k\ge 3$ and
$$
\hat J[\xi,\Phi]_{T,t}^{(2)}=\int\limits_t^T \xi_{t_2}d{\bf w}_{t_2}^{(i_2)}
\int\limits_{t_{2}}^{T}\Phi(t_1)d{\bf w}_{t_1}^{(i_1)}
\stackrel{\rm def}{=}
$$

\vspace{2mm}
$$
\stackrel{\rm def}{=}
\hbox{\vtop{\offinterlineskip\halign{
\hfil#\hfil\cr
{\rm l.i.m.}\cr
$\stackrel{}{{}_{N\to \infty}}$\cr
}} }\sum_{l=0}^{N-1}\xi_{\tau_l}\Delta{\bf w}_{\tau_{l}}^{(i_2)}
\int\limits_{\tau_{l+1}}^T
\Phi(t_1)d{\bf w}_{t_1}^{(i_1)}
$$

\vspace{5mm}
\noindent
for $k=2$. Here 
${\bf w}_{\tau}^{(i)}={\bf f}_{\tau}^{(i)}$ if $i=1,\ldots,m$ and
${\bf w}_\tau^{(0)}=\tau,$\ 
${\bf f}_{\tau}^{(i)}$ $(i=1,\ldots,m)$ are 
${\rm F}_{\tau}$-measurable for all $\tau\in[0,T]$ ($0\le t<T$)
independent standard Wiener processes,
$i_1,\ldots,i_k=0, 1,\ldots,m$.

Let us denote

\vspace{-1mm}
\begin{equation}
\label{2.18000000}
J[\xi,\Phi]_{T,t}^{(k)}=
\int\limits_t^T\ldots\int\limits_t^{t_{k-1}}
\Phi(t_1,\ldots,t_{k-1})\xi_{t_k}
d{\bf w}_{t_k}^{(i_k)}
\ldots d{\bf w}_{t_1}^{(i_1)},\ \ \ k\ge 2,
\end{equation}

\vspace{3mm}
\noindent
where the right-hand side of (\ref{2.18000000}) is the
iterated Ito stochastic 
integral.

Let us introduce the following iterated stochastic integrals

\vspace{2mm}
$$
\tilde J[\Phi]_{T,t}^{(k-1)}
=\int\limits_t^T d{\bf w}_{t_{k-1}}^{(i_{k-1})}\ldots 
\int\limits_{t_{3}}^T d{\bf w}_{t_{2}}^{(i_{2})}
\int\limits_{t_{2}}^{T}\Phi(t_1,t_2,\ldots,t_{k-1})
d{\bf w}_{t_1}^{(i_1)}
\stackrel{\rm def}{=}
$$

\vspace{2mm}
$$
\stackrel{\rm def}{=}\hbox{\vtop{\offinterlineskip\halign{
\hfil#\hfil\cr
{\rm l.i.m.}\cr
$\stackrel{}{{}_{N\to \infty}}$\cr
}} }\sum_{l=0}^{N-1}\Delta{\bf w}_{\tau_{l}}^{(i_{k-1})}
\int\limits_{\tau_{l+1}}^T
d{\bf w}_{t_{k-2}}^{(i_{k-2})}
\ldots\int\limits_{t_{3}}^T
d{\bf w}_{t_{2}}^{(i_{2})}
\int\limits_{t_{2}}^T
\Phi(t_1,t_2,\ldots,t_{k-1})d{\bf w}_{t_1}^{(i_1)},
$$

\vspace{4mm}
$$
J'[\Phi]_{T,t}^{(k-1)}=
\int\limits_t^T\ldots\int\limits_t^{t_{k-2}}
\Phi(t_1,\ldots,t_{k-1})
d{\bf w}_{t_{k-1}}^{(i_{k-1})}
\ldots d{\bf w}_{t_1}^{(i_1)},\ \ \ k\ge 2.
$$

\vspace{6mm}

Similarly to the proof of 
Theorem 1 it is easy to demonstrate that under the condition  
AII the stochastic  
integral 
$\tilde J[\Phi]_{T,t}^{(k-1)}$ exists and 

\begin{equation}
\label{432}
J'[\Phi]_{T,t}^{(k-1)}=\tilde J[\Phi]_{T,t}^{(k-1)}\ \ \ \hbox{w.~p.~1.}
\end{equation}

\vspace{3mm}

Moreover, using (\ref{432}) 
the following generalization
of Theorem 1 can be proved
similarly to the proof of Theorem 1.

\vspace{2mm}

{\bf Theorem 2}\ \cite{3}, \cite{3aaa} 
(also see \cite{77}-\cite{20}, \cite{12a}-\cite{12aaa}). 
{\it Suppose that the conditions  
{\rm AI, AII} of this section are fulfilled.
Then, the stochastic integral
$\hat J[\xi,\Phi]_{T,t}^{(k)}$ exists and for $k\ge 2$

\vspace{1mm}
$$
J[\xi,\Phi]_{T,t}^{(k)}=\hat J[\xi,\Phi]_{T,t}^{(k)}\ \ \ 
\hbox{w.~p.~{\rm 1.}}
$$
}

\vspace{2mm}

Let us consider the following stochastic integrals

\vspace{-1mm}
$$
I=\int\limits_t^{T}d{\bf f}_{t_2}^{(i_2)}\int\limits_{t_2}^T
\Phi_1(t_1,t_2)d{\bf f}_{t_1}^{(i_1)},\ \ \ 
J=\int\limits_t^{T}\int\limits_t^{t_2}
\Phi_2(t_1,t_2)d{\bf f}_{t_1}^{(i_1)}d{\bf f}_{t_2}^{(i_2)}.
$$

\vspace{2mm}

If we consider 
$$\int\limits_{t_2}^T
\Phi_1(t_1,t_2)d{\bf f}_{t_1}^{(i_1)}
$$ 

\vspace{2mm}
\noindent
as the integrand of $I$ 
and 

\vspace{-1mm}
$$
\int\limits_t^{t_2}
\Phi_2(t_1,t_2)d{\bf f}_{t_1}^{(i_1)}
$$ 

\vspace{2mm}
\noindent
as the integrand of $J$,
then, due to independence of these integrands we may mistakenly 
think that ${\sf M}\{IJ\}=0.$

But it is not the fact. Actually, using the 
integration order replacement technique in the stochastic
integral $I$, we have w.~p.~1 

$$
I=\int\limits_t^{T}\int\limits_t^{t_1}
\Phi_1(t_1,t_2)d{\bf f}_{t_2}^{(i_2)}d{\bf f}_{t_1}^{(i_1)}
=\int\limits_t^{T}\int\limits_t^{t_2}
\Phi_1(t_2,t_1)d{\bf f}_{t_1}^{(i_2)}d{\bf f}_{t_2}^{(i_1)}.
$$

\vspace{3mm}

So, using the standard properties of the Ito stochastic 
integral \cite{1}, we get

$$
{\sf M}\{IJ\}={\bf 1}_{\{i_1=i_2\}}\int\limits_t^{T}\int\limits_t^{t_2}
\Phi_1(t_2,t_1)\Phi_2(t_1,t_2)dt_1dt_2,
$$

\vspace{4mm}
\noindent
where ${\bf 1}_{\{A\}}$ is the indicator of the set $A$.

Let us consider the following statement.

\vspace{2mm}

{\bf Theorem 3}\ \cite{3}, \cite{3aaa} 
(also see \cite{77}-\cite{20}, \cite{12a}-\cite{12aaa}). {\it Let the 
conditions of Theorem {\rm 1} are fulfilled and $h(\tau)$  
is a continuous nonrandom function at the interval $[t,T]$.
Then

\begin{equation}
\label{2.19000000}
\int\limits_{t}^T \phi_\tau dw_\tau^{(k+1)}h(\tau) 
\hat I[\psi^{(k)}]_{T,\tau}=
\int\limits_{t}^T 
\phi_\tau h(\tau) dw_\tau^{(k+1)}\hat I[\psi^{(k)}]_{T,\tau}\ \ \
\hbox{{\rm w.~p.~1}},
\end{equation}

\vspace{3mm}
\noindent
where stochastic integrals on the left-hand side of {\rm (\ref{2.19000000})} 
as well as on the right-hand side of {\rm (\ref{2.19000000})} 
exist.} 

\vspace{2mm}

{\bf Proof.}\
According to Theorem 1, the iterated stochastic
integral on the right-hand side of (\ref{2.19000000}) exists.
In addition

\vspace{1mm}
$$
\int\limits_{t}^T \phi_\tau h(\tau) dw_\tau^{(k+1)} \hat I[\psi^{(k)}]_{T,\tau}=
\int\limits_{t}^T \phi_\tau dw_\tau^{(k+1)}h(\tau)\hat I[\psi^{(k)}]_{T,\tau} -
$$

$$
- \hbox{\vtop{\offinterlineskip\halign{
\hfil#\hfil\cr
{\rm l.i.m.}\cr
$\stackrel{}{{}_{N\to \infty}}$\cr
}} }\sum_{j=0}^{N-1} \phi_{\tau_j} \Delta h(\tau_{j})
\Delta w_{\tau_j}^{(k+1)} \hat I[\psi^{(k)}]_{T,\tau_{j+1}}\ \ \ 
\hbox{w. p. 1},
$$

\vspace{5mm}
\noindent
where $\Delta h(\tau_{j})=h(\tau_{j+1})-h(\tau_{j}).$

Using the arguments which resulted to the right equality
in (\ref{2.4000000}), we obtain

$$
\hbox{\vtop{\offinterlineskip\halign{
\hfil#\hfil\cr
{\rm l.i.m.}\cr
$\stackrel{}{{}_{N\to \infty}}$\cr
}} }\sum_{l=0}^{N-1}\phi_{\tau_l}\Delta h(\tau_{l})
\Delta w_{\tau_l}^{(k+1)} \hat I[\psi^{(k)}]_{T,\tau_{l+1}}=
$$

\begin{equation}
\label{2.20000001}
=\hbox{\vtop{\offinterlineskip\halign{
\hfil#\hfil\cr
{\rm l.i.m.}\cr
$\stackrel{}{{}_{N\to \infty}}$\cr
}} }G[\psi^{(k)}]_{N,0}
\sum_{l=0}^{j_k-1}\phi_{\tau_l}\Delta h(\tau_{l})
\Delta w_{\tau_l}^{(k+1)}\ \ \ \hbox{w. p. 1}.
\end{equation}

\vspace{3mm}

Using 
the Minkowski inequality, standard 
estimates for second moments of stochastic integrals as well as
continuity of the function $h(\tau)$, we obtain that
the second moment of the prelimit expression on  
the right-hand side of (\ref{2.20000001}) tends to zero 
when $N\to\infty.$ 
Theorem is proved.

Let us consider one corollary of Theorem 1.

\vspace{2mm}

{\bf Theorem 4}\ \cite{3}, \cite{3aaa} 
(also see \cite{77}-\cite{20}, \cite{12a}-\cite{12aaa}). {\it Under the conditions of Theorem {\rm 3}
the following equality is fulfilled                                

\vspace{-2mm}
$$
\int\limits_{t}^T h(t_1)\int\limits_{t}^{t_1}\phi_\tau dw_\tau^{(k+2)}
dw_{t_1}^{(k+1)} \hat I[\psi^{(k)}]_{T,t_1}=
$$

\begin{equation}
\label{2.21000000}
=\int\limits_{t}^T \phi_\tau dw_\tau^{(k+2)}\int\limits_{\tau}^T
h(t_1)dw_{t_1}^{(k+1)}\hat I[\psi^{(k)}]_{T,t_1}\ \ \ 
\hbox{w. p. {\rm 1}}.
\end{equation}

\vspace{4mm}

Moreover, the stochastic integrals in {\rm (\ref{2.21000000})} exist.}

\vspace{2mm}

{\bf Proof.}\ Using Theorem 1 two times,
we obtain

\vspace{1mm}
$$
\int\limits_{t}^T \phi_\tau dw_\tau^{(k+2)}\int\limits_{\tau}^T
h(t_1)dw_{t_1}^{(k+1)} \hat I[\psi^{(k)}]_{T,t_1}=
$$

\vspace{1mm}
$$
=\int\limits_{t}^T\psi_1(t_1)\ldots
\int\limits_{t}^{t_{k-1}}\psi_k(t_k)\int\limits_{t}^{t_k}
\rho_{\tau}dw_{\tau}^{(k+1)} dw_{t_k}^{(k)}\ldots dw_{t_1}^{(1)}=
$$

\vspace{1mm}
$$
=\int\limits_{t}^T \rho_\tau dw_\tau^{(k+1)}\int\limits_{\tau}^T\psi_k(t_k)
dw_{t_k}^{(k)}\ldots \int\limits_{t_{2}}^T\psi_1(t_1)dw_{t_1}^{(1)}\ \ \
\hbox{w. p. 1},
$$

\vspace{3mm}
\noindent
where
$$
\rho_{\tau}\stackrel{\rm def}{=}h(\tau)\int\limits_{t}^{\tau}
\phi_s dw_s^{(k+2)}.
$$ 

\vspace{2mm}

Theorem 4 is proved.

\vspace{5mm}

\section{Examples of Integration Order Replacement Technique
for the Concrete Iterated Ito Stochastic Integrals}

\vspace{5mm}

As we mentioned above, the formulas from this section
could be obtained using the Ito formula. However,
the method based on Theorem 1 is more simple
and familiar, since it deals with usual 
rules of the integration order replacement for Riemann integrals.

Using the integration order replacement technique
for iterated Ito stochastic integrals (Theorem 1), we obtain 
the following equalities which are fulfilled
w.~p.~1

\vspace{2mm}
$$
\int\limits_t^T\int\limits_t^{t_2}df_{t_1}dt_2=
\int\limits_t^T(T-t_1)df_{t_1},
$$

\vspace{2mm}
$$
\int\limits_t^T {\rm cos}(t_2-T)\int\limits_t^{t_2}df_{t_1}dt_2=\int\limits_t^T {\rm sin}(T-t_1)df_{t_1},
$$

\vspace{2mm}
$$
\int\limits_t^T {\rm sin}(t_2-T)\int\limits_t^{t_2}df_{t_1}dt_2=\int\limits_t^T 
\left({\rm cos}(T-t_1)-1\right)df_{t_1},
$$

\vspace{2mm}
$$
\int\limits_t^T e^{\alpha(t_2-T)}\int\limits_t^{t_2}df_{t_1}dt_2=
\frac{1}{\alpha}\int\limits_t^T\left(1-e^{\alpha(t_1-T)}\right)df_{t_1},\ \ \
\alpha\ne 0,
$$

\vspace{2mm}
$$
\int\limits_t^T(t_2-T)^{\alpha}\int\limits_t^{t_2}df_{t_1}dt_2=
-\frac{1}{\alpha+1}\int\limits_t^T(t_1-T)^{\alpha+1}df_{t_1},\ \ \ \alpha\ne -1,
$$

\vspace{2mm}
$$
J_{(100)T,t}=\frac{1}{2}\int\limits_t^T(T-t_1)^2 df_{t_1},
$$

\vspace{2mm}
$$
J_{(010)T,t}=\int\limits_t^T(t_1-t)(T-t_1)df_{t_1},
$$

\vspace{2mm}
\begin{equation}
\label{ex1}
J_{(110)T,t}=\int\limits_t^T(T-t_2)\int\limits_t^{t_2}df_{t_1}df_{t_2},\
\end{equation}

\vspace{2mm}
$$
J_{(101)T,t}=\int\limits_t^T\int\limits_t^{t_2}(t_2-t_1)df_{t_1}df_{t_2},
$$

\vspace{2mm}
$$
J_{(1011)T,t}=\int\limits_t^T\int\limits_t^{t_3}\int\limits_t^{t_2}(t_2-t_1)
df_{t_1}df_{t_2}df_{t_3},\
$$

\vspace{2mm}
$$
J_{(1101)T,t}=\int\limits_t^T\int\limits_t^{t_3}(t_3-t_2)\int\limits_t^{t_2}
df_{t_1}df_{t_2}df_{t_3},
$$

\vspace{2mm}
$$
J_{(1110)T,t}=\int\limits_t^T(T-t_3)\int\limits_t^{t_3}\int\limits_t^{t_2}
df_{t_1}df_{t_2}df_{t_3},
$$

\vspace{2mm}
$$
J_{(1100)T,t}=\frac{1}{2}\int\limits_t^T(T-t_2)^2\int\limits_t^{t_2}
df_{t_1}df_{t_2},
$$

\vspace{2mm}
$$
J_{(1001)T,t}=\frac{1}{2}\int\limits_t^T\int\limits_t^{t_2}(t_2-t_1)^2
df_{t_1}df_{t_2},\
$$

\vspace{2mm}
\begin{equation}
\label{ex2}
J_{(1010)T,t}=\int\limits_t^T(T-t_2)\int\limits_t^{t_2}(t_2-t_1)
df_{t_1}df_{t_2},
\end{equation}

\vspace{2mm}
$$
J_{(0110)T,t}=\int\limits_t^T(T-t_2)\int\limits_t^{t_2}(t_1-t)
df_{t_1}df_{t_2},
$$

\vspace{2mm}
$$
J_{(0101)T,t}=\int\limits_t^T\int\limits_t^{t_2}(t_2-t_1)(t_1-t)
df_{t_1}df_{t_2},
$$

\vspace{2mm}
$$
J_{(0010)T,t}=\frac{1}{2}\int\limits_t^T(T-t_1)(t_1-t)^2
df_{t_1},
$$

\vspace{2mm}
$$
J_{(0100)T,t}=\frac{1}{2}\int\limits_t^T(T-t_1)^2(t_1-t)
df_{t_1},
$$

\vspace{2mm}
$$
J_{(1000)T,t}=\frac{1}{3!}\int\limits_t^T(T-t_1)^3
df_{t_1},\
$$

\vspace{2mm}
$$
J_{(1\underbrace{{}_{0\ldots 0}}_{k-1})T,t}=
\frac{1}{(k-1)!}\int\limits_t^T(T-t_1)^{k-1}df_{t_1},
$$

\vspace{3mm}
$$
J_{(11\underbrace{{}_{0\ldots 0}}_{k-2})T,t}=
\frac{1}{(k-2)!}\int\limits_t^T(T-t_2)^{k-2}\int\limits_t^{t_2}df_{t_1}df_{t_2},
$$

\vspace{3mm}
$$
J_{(\underbrace{{}_{1\ldots 1}}_{k-1}0)T,t}=
\int\limits_t^T(T-t_1)J_{(\underbrace{{}_{1\ldots 1}}_{k-2})t_1,t}df_{t_1},
$$

\vspace{3mm}
$$
J_{(1\underbrace{{}_{0\ldots 0}}_{k-2}1)T,t}=
\frac{1}{(k-2)!}\int\limits_t^T\int\limits_t^{t_2}(t_2-t_1)^{k-2}df_{t_1}df_{t_2},
$$

\vspace{3mm}
$$
J_{(10\underbrace{{}_{1\ldots 1}}_{k-2})T,t}=
\int\limits_t^T\ldots \int\limits_t^{t_3}\int\limits_t^{t_2}(t_2-t_1)
df_{t_1}df_{t_2}\ldots df_{t_{k-1}},
$$

\vspace{5mm}
$$
J_{(\underbrace{{}_{1\ldots 1}}_{k-2}01)T,t}=
\int\limits_t^T\int\limits_t^{t_{k-1}}(t_{k-1}-t_{k-2})
\int\limits_t^{t_{k-2}}\ldots
\int\limits_t^{t_2}df_{t_1}\ldots df_{t_{k-3}}df_{t_{k-2}}df_{t_{k-1}},
$$

\vspace{5mm}
$$
J_{(10)T,t}+J_{(01)T,t}=(T-t)J_{(1)T,t},\
$$
        
\vspace{2mm}
$$
J_{(110)T,t}+J_{(101)T,t}+J_{(011)T,t}=(T-t)J_{(11)T,t},\
$$

\vspace{2mm}
$$
J_{(001)T,t}+J_{(010)T,t}+J_{(100)T,t}=\frac{(T-t)^2}{2}
J_{(1)T,t},
$$

\vspace{6mm}
$$
J_{(1100)T,t}+J_{(1010)T,t}+J_{(1001)T,t}+
J_{(0110)T,t}+
$$

\vspace{2mm}
$$
+J_{(0101)T,t}+J_{(0011)T,t}=\frac{(T-t)^2}{2}
J_{(11)T,t},
$$

\vspace{6mm}
$$
J_{(1000)T,t}+J_{(0100)T,t}+J_{(0010)T,t}+J_{(0001)T,t}=\frac{(T-t)^3}{3!}
J_{(1)T,t},
$$

\vspace{3mm}
$$
J_{(1110)T,t}+J_{(1101)T,t}+J_{(1011)T,t}+J_{(0111)T,t}=
(T-t)J_{(111)T,t},
$$

\vspace{4mm}
$$
\sum_{l=1}^k J_{(\underbrace{{}_{0\ldots 0}}_{l-1}1 
\underbrace{{}_{0\ldots 0}}_{k-l})T,t}=\frac{1}{(k-1)!}(T-t)^{k-1}J_{(1)T,t},
$$

\vspace{4mm}
$$
\sum_{l=1}^k J_{(\underbrace{{}_{1\ldots 1}}_{l-1}0 
\underbrace{{}_{1\ldots 1}}_{k-l})T,t}=(T-t)
J_{(\underbrace{{}_{1\ldots 1}}_{k-1} )T,t},
$$

\vspace{4mm}
$$
\sum_{{}_{\stackrel{l_1+\ldots+l_k=m}{l_i\in\{0,\ 1\},\ i=1,\ldots,k}}}
J_{(l_1\ldots l_k)T,t}=\frac{(T-t)^{k-m}}{(k-m)!}
J_{(\underbrace{{}_{1\ldots 1}}_{m})T,t},
$$

\vspace{4mm}
\noindent
where
$$
J_{(l_1\ldots l_k)T,t}=\int\limits_t^T\ldots \int\limits_t^{t_2}dw_{t_1}
^{(1)}\ldots dw_{t_k}^{(k)},
$$

\vspace{3mm}
\noindent
$l_i=1$ when $w_{t_i}^{(i)}=f_{t_i}$ and $l_i=0$ when
$w_{t_i}^{(i)}=t_i$ $(i=1,\ldots,k),$
$f_\tau$ is a standard Wiener process.

Let us consider two examples and show explicitly the technique on 
integration
order replacement for iterated Ito stochastic integrals.

\vspace{2mm}

{\bf Example~1.} {\it Let us prove the equality {\rm (\ref{ex1}).} 
Using Theorems {\rm 1} and {\rm 3,} 
we obtain}

\vspace{2mm}
$$
J_{(110)T,t}\stackrel{\rm def}{=}
\int\limits_t^T\int\limits_t^{t_3}\int\limits_t^{t_2} df_{t_1}df_{t_2}dt_3=
$$

\vspace{2mm}
$$
=\int\limits_t^T df_{t_1}\int\limits_{t_1}^T df_{t_2} \int\limits_{t_2}^T dt_3=
$$

\vspace{2mm}
$$
=
\int\limits_t^T df_{t_1}\int\limits_{t_1}^T df_{t_2}(T-t_2)=
$$

\vspace{2mm}
$$
=\int\limits_t^T df_{t_1}\int\limits_{t_1}^T (T-t_2) df_{t_2}=
$$

\vspace{2mm}
$$
=
\int\limits_t^T(T-t_2)\int\limits_t^{t_2}df_{t_1}df_{t_2}\ \ \ 
\hbox{w.~p.~1.}
$$
 
\vspace{5mm}

{\bf Example~2.} {\it Let us prove the equality {\rm (\ref{ex2}).} 
Using Theorems {\rm 1} and {\rm 3,} 
we obtain}

\vspace{2mm}
$$
J_{(1010)T,t}\stackrel{\rm def}{=}
\int\limits_t^T\int\limits_t^{t_4}
\int\limits_t^{t_3}\int\limits_t^{t_2}df_{t_1}dt_2df_{t_3}dt_4=
$$

\vspace{2mm}
$$
=
\int\limits_t^T df_{t_1}\int\limits_{t_1}^T dt_2 \int\limits_{t_2}^T df_{t_3}
\int\limits_{t_3}^T dt_4
=
$$

\vspace{2mm}
$$
=\int\limits_t^T df_{t_1}\int\limits_{t_1}^T 
dt_2 \int\limits_{t_2}^T df_{t_3}(T-t_3)=
$$

\vspace{2mm}
$$
=\int\limits_t^T df_{t_1}
\int\limits_{t_1}^T dt_2 \int\limits_{t_2}^T(T-t_3)df_{t_3}=
$$

\vspace{2mm}
$$
=
\int\limits_t^T(T-t_3)
\int\limits_t^{t_3}\int\limits_t^{t_2} df_{t_1}dt_2df_{t_3}=
$$

\vspace{2mm}
$$
=\int\limits_t^T(T-t_3)\left(\int\limits_t^{t_3}
\int\limits_t^{t_2} df_{t_1}dt_2\right)df_{t_3}=
$$

\vspace{2mm}
$$
=
\int\limits_t^T(T-t_3)\left(
\int\limits_t^{t_3}df_{t_1}\int\limits_{t_1}^{t_3}dt_2\right)df_{t_3}=
$$

\vspace{2mm}
$$
=\int\limits_t^T(T-t_3)\left(
\int\limits_t^{t_3}df_{t_1}(t_3-t_1)\right)df_{t_3}=
$$

\vspace{2mm}
$$
=
\int\limits_t^T(T-t_3)\left(
\int\limits_t^{t_3}(t_3-t_1)df_{t_1}\right)df_{t_3}=
$$

\vspace{2mm}
$$
=\int\limits_t^T(T-t_2)\int\limits_t^{t_2}(t_2-t_1)
df_{t_1}df_{t_2}\ \ \ \hbox{w. p. 1.}
$$

\vspace{5mm}

\section{Integration Order Replacement Technique for Iterated
Sto\-chas\-tic Integrals With Respect to Martingale}

\vspace{5mm}

In this section, we will generalize the theorems on
integration order replacement for iterated Ito stochastic 
integrals to the class of iterated stochastic 
integrals with respect to martingale.

Let
$(\Omega,{\rm F},{\sf P})$ be a complete probability space 
and 
let $\{{\rm F}_t, t\in[0, T]\}$ be a nondecreasing 
family of $\sigma$-algebras defined on the probability space
$(\Omega,{\rm F},{\sf P}).$
Suppose that $M_t,$ $t\in[0,T]$ is an ${\rm F}_t$-measurable 
martingale for all 
$t\in[0,T]$, which 
satisfies
the condition
${\sf M}\left\{\left|M_t\right|\right\}<\infty$. Moreover,
for all $t\in[0,T]$ there exists
an ${\rm F}_t$- measurable and 
nonnegative w.~p.~1 stochastic process 
$\rho_t,$ $t\in[0,T]$ such that 

\vspace{-1mm}
$$
{\sf M}\left\{\left(M_s-M_t\right)^2\left|\right.{\rm F}_t\right\}=
{\sf M}\left\{\int\limits_t^s\rho_{\tau}d\tau \biggl|\biggr.
{\rm F}_t\right\}\ \ \
\hbox{w.\ p.\ 1},
$$

\vspace{2mm}
\noindent
where $0\le t<s\le T.$

Let us consider the class $H_2(\rho,[0,T])$
of stochastic processes
$\varphi_t,$ $t\in[0,T]$, which are
${\rm F}_t$-measurable for all 
$t\in[0,T]$ and 
satisfy
the condition 

\vspace{-1mm}
$$
{\sf M}\left\{\int\limits_0^T\varphi_t^2 \rho_t dt\right\}<\infty.
$$

\vspace{2mm}

For any partition 
$\tau_j^{(N)},$ $j=0, 1, \ldots, N$ of
the interval $[0,T]$ such that

\begin{equation}
\label{w11ggg}
0=\tau_0^{(N)}<\tau_1^{(N)}<\ldots <\tau_N^{(N)}=T,\ \ \ \
\max\limits_{0\le j\le N-1}\left|\tau_{j+1}^{(N)}-\tau_j^{(N)}\right|\to 0\ \
\hbox{if}\ \ N\to \infty
\end{equation}

\vspace{3mm}
\noindent
we will define the sequence of step functions 

$$
\varphi^{(N)}(t,\omega)=
\varphi_j(\omega)\ \ \ \hbox{w. p. 1}\ \ \ 
\hbox{for}\ \ \ t\in\left[\tau_j^{(N)},\tau_{j+1}^{(N)}\right),
$$

\vspace{3mm}
\noindent
where $j=0, 1,\ldots,N-1,$ $N=1, 2,\ldots$

Let us define the stochastic integral with respect to martingale
for
$\varphi(t,\omega)\in H_2(\rho,[0,T])$ 
as the 
following mean-square limit \cite{1}

\vspace{-1mm}
$$
\hbox{\vtop{\offinterlineskip\halign{
\hfil#\hfil\cr
{\rm l.i.m.}\cr
$\stackrel{}{{}_{N\to \infty}}$\cr
}} }\sum_{j=0}^{N-1}\varphi^{(N)}\left(\tau_j^{(N)},\omega\right)
\left(M\left(\tau_{j+1}^{(N)},\omega\right)-
M\left(\tau_j^{(N)},\omega\right)\right)
\stackrel{\rm def}{=}\int\limits_0^T\varphi_\tau dM_\tau,
$$

\vspace{2mm}
\noindent
where $\varphi^{(N)}(t,\omega)$ is any step function
from the class $H_2(\rho,[0,T])$,
which converges
to the function $\varphi(t,\omega)$
in the following sense

\vspace{-1mm}
$$
\hbox{\vtop{\offinterlineskip\halign{
\hfil#\hfil\cr
{\rm lim}\cr
$\stackrel{}{{}_{N\to \infty}}$\cr
}} }\int\limits_0^T{\sf M}\left\{\left|
\varphi^{(N)}(t,\omega)-\varphi(t,\omega)\right|^2\right\}\rho_tdt=0.
$$

\vspace{2mm}

It is well known  \cite{1} that the stochastic integral

\vspace{-1mm}
$$
\int\limits_0^T\varphi_{\tau} dM_{\tau}
$$

\vspace{2mm}
\noindent
exists and it does not depend on the selection 
of sequence 
$\varphi^{(N)}(t,\omega)$.

Let $\tilde H_2(\rho,[0,T])$ be the class of 
stochastic processes  $\varphi_{\tau},$ $\tau\in[0,T],$
which are
mean-square
continuous for all $\tau\in[0,T]$ and belong to the 
class $H_2(\rho,[0,T])$.

Let us consider the following iterated stochastic integrals

\vspace{-1mm}
\begin{equation}
\label{32.001}
S[\phi,\psi^{(k)}]_{T,t}=\int\limits_{t}^{T}\psi_1(t_1)\ldots
\int\limits_t^{t_{k-1}}\psi_k(t_k)\int\limits_t^{t_k}
\phi_{\tau}dM^{(k+1)}_{\tau}dM^{(k)}_{t_k}
\ldots dM^{(1)}_{t_1},
\end{equation}

\begin{equation}
\label{32.002}
S[\psi^{(k)}]_{T,t}=\int\limits_{t}^{T}\psi_1(t_1)\ldots
\int\limits_t^{t_{k-1}}\psi_k(t_k)dM^{(k)}
_{t_k}
\ldots dM^{(1)}_{t_1}.
\end{equation}

\vspace{2mm}
\noindent
Here $\phi_\tau\in \tilde H_2(\rho,[t,T])$ and
$\psi_1(\tau),\ldots,\psi_k(\tau)$ are continuous nonrandom 
functions at the interval
$[t,T]$,
$M_\tau^{(l)}=M_\tau$ or $M_\tau^{(l)}=\tau$
if $\tau\in[t,T],$
$l=1,\ldots,k+1,$\ $M_{\tau}$ is the martingale defined above.

Let us define the iterated stochastic 
integral 
$\hat S[\psi^{(k)}]_{T,s},$ $0\le t\le s\le T,$ $k\ge 1$
with respect to martingale
$$
\hat S[\psi^{(k)}]_{T,s}=\int\limits_s^T\psi_k(t_k)dM_{t_k}^{(k)}
\ldots \int\limits_{t_{2}}^T \psi_1(t_1)dM_{t_1}^{(1)}
$$

\vspace{3mm}
\noindent
by the following recurrence 
relation

\vspace{-1mm}
\begin{equation}
\label{32.003}
\hat S[\psi^{(k)}]_{T,t}
\stackrel{\rm def}{=}
\hbox{\vtop{\offinterlineskip\halign{
\hfil#\hfil\cr
{\rm l.i.m.}\cr
$\stackrel{}{{}_{N\to \infty}}$\cr
}} }
\sum^{N-1}_{l=0} \psi_k(\tau_{l})\Delta M_{\tau_l}^{(k)} 
\hat S[\psi^{(k-1)}]_{T,\tau_{l+1}},
\end{equation}

\vspace{3mm}
\noindent
where $k\ge 1,$ 
$\hat S[\psi^{(0)}]_{T,s}\stackrel{\rm def}{=}1,$
$[s,T]\subseteq[t,T],$\ here and further $\Delta M_{\tau_l}^{(i)}=
M_{\tau_{l+1}}^{(i)}-M_{\tau_l}^{(i)},$
$i=1,\ldots,k+1,$ $l=0, 1,\ldots,N-1,$
$\{\tau_l\}_{l=0}^N$ is the partition of the interval $[t, T]$,
which satisfies the condition similar to (\ref{w11ggg}),
another notations are the same as in
(\ref{32.001}), (\ref{32.002}).

Further,
let us define the iterated stochastic 
integral 
$\hat S[\phi,\psi^{(k)}]_{T,t},$ $k\ge 1$ of the form

$$
\hat S[\phi,\psi^{(k)}]_{T,t}=\int\limits_{t}^T \phi_s dM_s^{(k+1)}
\hat S[\psi^{(k)}]_{T,s}
$$

\vspace{2mm}
\noindent
by the equality

$$
\hat S[\phi,\psi^{(k)}]_{T,t}
\stackrel{\rm def}{=}
\hbox{\vtop{\offinterlineskip\halign{
\hfil#\hfil\cr
{\rm l.i.m.}\cr
$\stackrel{}{{}_{N\to \infty}}$\cr
}} }
\sum^{N-1}_{l=0} \phi_{\tau_{l}}\Delta M_{\tau_l}^{(k+1)} 
\hat S[\psi^{(k)}]_{T,\tau_{l+1}},
$$

\vspace{4mm}
\noindent
where the sense of notations included in 
(\ref{32.001})--(\ref{32.003}) is saved.

Let us formulate the theorem on integration order replacement 
for the iterated
stochastic integrals with respect to martingale,
which is the generalization of Theorem 1.

\vspace{2mm}

{\bf Theorem 5}\ \cite{3}, \cite{3aaa} 
(also see \cite{77}-\cite{20}, \cite{12a}-\cite{12aaa}). {\it Let 
$\phi_\tau\in\tilde H_2(\rho,[t,T]),$ every $\psi_l(\tau)$
$(l=1,\ldots,k)$ is a continuous nonrandom function at the interval 
$[t,T],$ and
$|\rho_{\tau}|\le K<\infty $ w.~p.~{\rm 1} for all
$\tau\in[t,T].$
Then, the stochastic integral 
$\hat S[\phi,\psi^{(k)}]_{T,t}$ exists and

\vspace{1mm}
$$
S[\phi,\psi^{(k)}]_{T,t}=\hat S[\phi,\psi^{(k)}]_{T,t}\ \ \ 
\hbox{w.~p.~{\rm 1}}.
$$ 
}

\vspace{1mm}

The proof of Theorem 5 is similar to the proof of 
Theorem 1.

\vspace{2mm}

{\bf Remark 2.}\ {\it Let us note that we can propose 
another variant of the conditions in Theorem {\rm 5}. For example,
if we not require the boundedness of the process  
$\rho_{\tau}$, then it is necessary to require
the fulfillment of the following additional 
conditions{\rm :}

\vspace{2mm}

{\rm 1}. ${\sf M}\{|\rho_{\tau}|\}<\infty$ for all $\tau\in[t,T].$

\vspace{2mm}

{\rm 2}. The process $\rho_{\tau}$ is independent with 
the processes  $\phi_{\tau}$
and $M_{\tau}.$
}

\vspace{2mm}

{\bf Remark 3.}\ {\it Note that it is well known
the construction of stochastic integral
with respect to the Wiener process with integrable process, which is
not an ${\rm F}_{\tau}$-measurable 
stochastic  
process  --- the so-called 
Stratonovich stochastic integral {\rm \cite{2}}.

The stochastic integral $\hat S[\phi,\psi^{(k)}]_{T,t}$
is also the stochastic  
integral with 
integrable process, which is 
not an ${\rm F}_{\tau}$-measurable 
stochastic  
process. 
However, under the conditions of Theorem {\rm 5}

$$
S[\phi,\psi^{(k)}]_{T,t}=\hat S[\phi,\psi^{(k)}]_{T,t}\ \ \ 
\hbox{w. p. {\rm 1,}}
$$

\vspace{3mm}
\noindent
where $S[\phi,\psi^{(k)}]_{T,t}$ is a usual
iterated stochastic integral with respect to martingale.
If, for example, $M_{\tau}, \tau\in[t,T]$
is the Wiener process, then the question on connection between stochastic 
integral $\hat S[\phi,\psi^{(k)}]_{T,t}$
and Stratonovich stochastic integral
is solving as a standard question on connection
between Stratonovich and Ito stochastic integrals {\rm \cite{2}}.
}  

\vspace{2mm}

Let us consider several statements, which are the generalizations of theorems
formulated in the previous sections.

Assume that $D_k=\{(t_1,\ldots,t_k):\ t\le t_1<\ldots<t_k\le T\}$
and the following conditions are met:

\vspace{2mm}

BI.\ $\xi_{\tau}\in \tilde H_2(\rho,[t,T]).$

\vspace{2mm}

BII.\ $\Phi(t_1,\ldots,{t_{k-1}})$ is a continuous nonrandom function
in the closed
domain $D_{k-1}$ (recall that 
we use the same symbol $D_{k-1}$ to denote the open and closed 
domains corresponding to the domain $D_{k-1};$
however, we always specify what domain we consider (open or closed)).

\vspace{2mm}

Let us define the following stochastic integrals
with respect to martingale

\vspace{1mm}
$$
\hat S[\xi,\Phi]_{T,t}^{(k)}=\int\limits_t^T \xi_{t_k}dM_{t_k}^{(k)}\ldots 
\int\limits_{t_{3}}^T dM_{t_{2}}^{(2)}
\int\limits_{t_{2}}^{T}\Phi(t_1,t_2,\ldots,t_{k-1})
dM_{t_1}^{(1)}
\stackrel{\rm def}{=}
$$

\vspace{3mm}
$$
\stackrel{\rm def}{=}\hbox{\vtop{\offinterlineskip\halign{
\hfil#\hfil\cr
{\rm l.i.m.}\cr
$\stackrel{}{{}_{N\to \infty}}$\cr
}} }\sum_{l=0}^{N-1}\xi_{\tau_l}\Delta M_{\tau_{l}}^{(k)}
\int\limits_{\tau_{l+1}}^T
dM_{t_{k-1}}^{(k-1)}
\ldots\int\limits_{t_{3}}^T
dM_{t_{2}}^{(2)}
\int\limits_{t_{2}}^T
\Phi(t_1,t_2,\ldots,t_{k-1})dM_{t_1}^{(1)}
$$

\vspace{5mm}
\noindent
for $k\ge 3$ and
$$
\hat S[\xi,\Phi]_{T,t}^{(2)}=\int\limits_t^T \xi_{t_2}dM_{t_2}^{(2)}
\int\limits_{t_{2}}^{T}\Phi(t_1)dM_{t_1}^{(1)}
\stackrel{\rm def}{=}
$$

\vspace{2mm}
$$
\stackrel{\rm def}{=}
\hbox{\vtop{\offinterlineskip\halign{
\hfil#\hfil\cr
{\rm l.i.m.}\cr
$\stackrel{}{{}_{N\to \infty}}$\cr
}} }\sum_{l=0}^{N-1}\xi_{\tau_l}\Delta M_{\tau_{l}}^{(2)}
\int\limits_{\tau_{l+1}}^T
\Phi(t_1)dM_{t_1}^{(1)}
$$

\vspace{6mm}
\noindent
for $k=2$, where the sense of notations included 
in (\ref{32.001})--(\ref{32.003}) is saved.
Moreover,
the stochastic process  $\xi_{\tau},$ $\tau\in[t,T]$
belongs to the class $\tilde H_2(\rho,[t,T]).$

In addition, let

\vspace{-1mm}
\begin{equation}
\label{32.006}
S[\xi,\Phi]_{T,t}^{(k)}=
\int\limits_t^T\ldots\int\limits_t^{t_{k-1}}
\Phi(t_1,\ldots,t_{k-1})\xi_{t_k}
dM_{t_k}^{(k)}
\ldots dM_{t_1}^{(1)},\ \ \ k\ge 2,
\end{equation}

\vspace{3mm}
\noindent
where the right-hand side of (\ref{32.006}) is the iterated stochastic 
integral with respect to martingale.

Let us introduce the following iterated stochastic integrals
with respect to martingale

$$
\tilde S[\Phi]_{T,t}^{(k-1)}
=\int\limits_t^T d{M}_{t_{k-1}}^{(k-1)}
\ldots 
\int\limits_{t_{3}}^T d{M}_{t_{2}}^{(2)}
\int\limits_{t_{2}}^{T}\Phi(t_1,t_2,\ldots,t_{k-1})
d{M}_{t_1}^{(1)}
\stackrel{\rm def}{=}
$$

\vspace{2mm}
$$
\stackrel{\rm def}{=}\hbox{\vtop{\offinterlineskip\halign{
\hfil#\hfil\cr
{\rm l.i.m.}\cr
$\stackrel{}{{}_{N\to \infty}}$\cr
}} }\sum_{l=0}^{N-1}\Delta{M}_{\tau_{l}}^{(k-1)}
\int\limits_{\tau_{l+1}}^T
d{M}_{t_{k-2}}^{(k-2)}
\ldots\int\limits_{t_{3}}^T
d{M}_{t_{2}}^{(2)}
\int\limits_{t_{2}}^T
\Phi(t_1,t_2,\ldots,t_{k-1})d{M}_{t_1}^{(1)},
$$

\vspace{2mm}
$$
S'[\Phi]_{T,t}^{(k-1)}=
\int\limits_t^T\ldots\int\limits_t^{t_{k-2}}
\Phi(t_1,\ldots,t_{k-1})
d{M}_{t_{k-1}}^{(k-1)}
\ldots d{M}_{t_1}^{(1)},\ \ \ k\ge 2.
$$

\vspace{4mm}

It is easy to demonstrate similarly to the proof of 
Theorem 5 that 
under the condition  
BII the stochastic  
integral
$\tilde S[\Phi]_{T,t}^{(k-1)}$ exists
and

$$
S'[\Phi]_{T,t}^{(k-1)}=\tilde S[\Phi]_{T,t}^{(k-1)}\ \ \ 
\hbox{w.~p.~{\rm 1}}.
$$ 

\vspace{3mm}

In its turn, using this fact we can
prove the following theorem
similarly to the proof of Theorem 5.

\vspace{2mm}

{\bf Theorem 6}\ \cite{3}, \cite{3aaa} 
(also see \cite{77}-\cite{20}, \cite{12a}-\cite{12aaa}).
{\it Let the conditions {\rm BI, BII} of this section are fulfilled
and $|\rho_{\tau}|\le K<\infty $ w.~p.~{\rm 1} for all 
$\tau\in[t,T].$
Then, the stochastic integral
$\hat S[\xi,\Phi]_{T,t}^{(k)}$ exists and for $k\ge 2$

$$
S[\xi,\Phi]_{T,t}^{(k)}=\hat S[\xi,\Phi]_{T,t}^{(k)}\ \ \ 
\hbox{w.~p.~{\rm 1}}.
$$
}

\vspace{2mm}

Theorem 6 is the generalization of Theorem 2 for the case  
of iterated stochastic integrals with respect to martingale.

Let us  consider two statements. 

\vspace{2mm}
 
{\bf Theorem 7}\ \cite{3}, \cite{3aaa} 
(also see \cite{77}-\cite{20}, \cite{12a}-\cite{12aaa}). {\it
Let the conditions of Theorem {\rm 5} are fulfilled and $h(\tau)$ 
is a continuous nonrandom function at the interval $[t,T]$.
Then

\begin{equation}
\label{32.007}
\int\limits_{t}^T \phi_\tau dM_\tau^{(k+1)}h(\tau) 
\hat S[\psi^{(k)}]_{T,\tau}=
\int\limits_{t}^T \phi_\tau h(\tau) dM_\tau^{(k+1)}
\hat S[\psi^{(k)}]_{T,\tau}\ \ \ \hbox{w. p. {\rm 1}}
\end{equation}

\vspace{3mm}
\noindent
and stochastic  
integrals on the left-hand side of {\rm (\ref{32.007})}
as well as on the right-hand side of {\rm (\ref{32.007})}
exist.} 

\vspace{2mm}

{\bf Theorem 8}\ \cite{3}, \cite{3aaa} 
(also see \cite{77}-\cite{20}, \cite{12a}-\cite{12aaa}). {\it Under the 
conditions of Theorem {\rm 7}

$$
\int\limits_{t}^T h(t_1)\int\limits_{t}^{t_1}\phi_\tau dM_\tau^{(k+2)}
dM_{t_1}^{(k+1)} \hat S[\psi^{(k)}]_{T,t_1}=
$$

\begin{equation}
\label{32.008}
=\int\limits_{t}^T \phi_\tau dM_\tau^{(k+2)}\int\limits_{\tau}^T
h(t_1)dM_{t_1}^{(k+1)}\hat S[\psi^{(k)}]_{T,t_1}\ \hbox{\rm w.\ p.\ 1}.
\end{equation}

\vspace{4mm}

Moreover, the stochastic integrals in {\rm (\ref{32.008})} exist.} 

\vspace{2mm}

The proofs of Theorems 7 and 8 are similar to the proofs of 
Theorems 3 and 4 correspondingly.

\vspace{2mm}

{\bf Remark 4.}\
{\it The integration order replacement technique for iterated Ito stochastic 
integrals {\rm (}Theorems {\rm 1--4)} \cite{77}-\cite{50}, \cite{12a}-\cite{12aaa}
has been successfully applied for construction of the so-called
unified Taylor--Ito and Taylor--Stratonovich expansions \cite{12a}-\cite{12aaa}
{\rm (}see references therein{\rm )}
as well as
for proof and development of the mean-square 
approximation method for iterated 
Ito and Stratonovich stochastic integrals
based on generalized multiple Fourier series \cite{12a}-\cite{12aaa}
{\rm (}see references therein{\rm )}.}

\end{document}